\documentclass[a4paper,10pt]{elsarticle}
\usepackage[utf8]{inputenc}

\usepackage{a4wide}
\usepackage{amsmath,amssymb}
\usepackage{enumerate}
\usepackage{color}
\usepackage{bbold}

\usepackage{graphicx}
\usepackage{caption}
\usepackage{subcaption}

\usepackage{lineno,hyperref}

\newtheorem{theorem}{Theorem}
\newtheorem{lemma}{Lemma}
\newtheorem{definition}{Definition}
\newtheorem{proposition}{Proposition}
\newtheorem{remark}{Remark}
\newtheorem{corollary}{Corollary}

\newcommand{\proof}{\noindent {\bf Proof. }}
\newcommand{\eproof}{\hfill $\square$}

\newcommand{\uu}{\mathbf{u}}

\newcommand{\ww}{\mathbf{w}}
\newcommand{\bd}{\boldsymbol{d}}

\newcommand{\bX}{\boldsymbol{X}}
\newcommand{\bF}{\boldsymbol{F}}
\newcommand{\bG}{\boldsymbol{G}}
\newcommand{\bQ}{\boldsymbol{Q}}

\newcommand{\bXi}{\boldsymbol{\Xi}}

\newcommand{\Sf}{\boldsymbol{S}}

\newcommand{\Pib}{\mathfrak{X}}

\begin{document}

\begin{frontmatter}

\title{Asymptotic ensemble stabilizability of the Bloch equation}


\author[mymainaddress]{Francesca C. Chittaro}
\ead{francesca.chittaro@univ-tln.fr}
\author[mymainaddress]{Jean-Paul Gauthier}
\ead{gauthier@univ-tln.fr}

\address[mymainaddress]{Aix Marseille Universit\'{e}, CNRS, ENSAM, LSIS UMR 7296, 13397 Marseille, France, and Universit\'e de Toulon,
CNRS, LSIS UMR 7296, 83957 La Garde, France}

\begin{abstract}
In this paper we are concerned with  the stabilizability to an equilibrium point of an ensemble of non
interacting half-spins. We assume that the spins are immersed in a static magnetic field,  with dispersion in the Larmor frequency, and are controlled by a time varying 
transverse field. Our goal is to steer the whole ensemble to the uniform 
``down'' position.

Two cases are addressed:
for a finite ensemble of spins, we provide a control function (in feedback form)
that asymptotically stabilizes  the ensemble to the ``down'' position, generically with respect to the initial condition.
For an ensemble containing a countable number of spins, we construct  a sequence of control functions such that the sequence of the corresponding 
solutions pointwise converges, asymptotically
in time, to the target state, generically with respect to the  initial conditions.

The control functions proposed are uniformly bounded and continuous.
\end{abstract}


\end{frontmatter}

\section{Introduction}

Ensemble controllability (also called \emph{simultaneous controllability}) is a notion introduced in \cite{Li-Khaneja2006,Li-Khanejacdc2006,Li-Khaneja2009} for quantum systems described by a family of parameter-dependent
ordinary 
differential equations; 
it concerns the possibility of finding control functions that compensate the dispersion in the parameters and drive the whole family (ensemble) from some initial state to some prescribed target state.

Such an issue is motivated by recent engineering applications, such as, for instance,
quantum control (see for instance \cite{Li-Khaneja2009,beauchard-coron-rouchon2010,leghtas-sarlette-rouchon,altafini-ens-scl} and references therein),  distributed parameters 
systems and PDEs 
\cite{dahleh,chambrioncdc2013,coron-book,curtain,helmke-schoenlein2014}, and flocks of identical systems \cite{Brockett2010}. 

General results for the ensemble controllability of linear  and nonlinear systems, in continuous and discrete time, can be found in the recent papers
\cite{Li2011,dirr,agra-ensemble,turi-ens,helmke-schoenlein2015}.

This paper deals with the simultaneous control of an ensemble of half-spins
immersed on a magnetic field, where each spin is described by a magnetization vector $\boldsymbol{M}\in \mathbb{R}^3$, subject to 
the dynamics
$\frac{d\boldsymbol{M}}{dt}=-\gamma \boldsymbol{M}\times \boldsymbol{B}(\boldsymbol{r},t)$,  
where $\boldsymbol{B}(\boldsymbol{r},t)$ is a magnetic field composed by a static component directed along the $z$-axis, and 
a time varying component on the $xy$-plane, called \emph{radio-frequency (rf) field}, and $\gamma$ denotes the gyromagnetic ratios of the spins.
In this system, since all spins are controlled by the same magnetic 
field $\boldsymbol{B}(\boldsymbol{r},t)$, the spatial dispersion in the amplitude of the magnetic field gives rise to the following
 inhomogeneities in the dynamics:  \emph{rf inhomogeneity}, caused by dispersion in the radio-frequency field, and a spread in the Larmor frequency,
 given by dispersion of the static component of the field.
This problem arises, for instance, in NMR spectroscopy (see \cite{Glaser2015} and  references in \cite{altafini2007,Li-Khaneja2009,beauchard-coron-rouchon2010}).

The task of controlling such system is wide, multi-faceted and very rich, depending on the cardinality of the set of the spin to be controlled (and the topology of this set), on the 
particular notion of controllability addressed, and on the functional space where control functions live.

The above-cited articles \cite{Li-Khaneja2006,Li-Khanejacdc2006,Li-Khaneja2009} are concerned with both rf inhomogeneity and Larmor dispersion, with dispersion parameters that belong to some compact domain $\mathcal{D}$. 
The magnetization vector
of the system is thus a function on $\mathcal{D}$, taking values in the unit sphere of $\mathbb{R}^3$, and ensemble controllability has to be intended as convergence in the 
$L^{\infty}(\mathcal{D},\mathbb{R}^3)$-norm.  
The controllability result is achieved by means of Lie algebraic techniques coupled with adiabatic evolution, and holds for both bounded and unbounded controls.

In \cite{beauchard-coron-rouchon2010}, the authors focus on systems subject to Larmor dispersions, and provide a complete analysis of controllability properties of the ensemble in
different scenarios, such as: bounded/unbounded controls; finite time/asymptotic controllability; approximate/exact controllability in the $L^{2}(\mathcal{D},\mathbb{R}^3)$ norm; boundedness/unboundedness of the set 
$\mathcal{D}$. In particular, results on exact local controllability with unbounded controls are provided.

In this paper we consider an ensemble of Bloch equations presenting Larmor dispersion, with frequencies belonging to some bounded subset $\mathcal{E}\subset \mathbb{R}$.  
Coupling a Lyapunov function approach with some tools of dynamical systems theory, we exhibit a control function (in feedback form) that approximately
drives, asymptotically in time and generically with respect to the initial conditions, all spins to the ``down'' position.
Two cases are addressed: if the set $\mathcal{E}$ is finite, our strategy provides exact exponential stabilizability in infinite time, while in the case where $\mathcal{E}$ is a countable 
collection of energies, our approach implies asymptotic pointwise convergence towards the target state.

Feedback control is a widely used tool for stabilization of control-affine systems (see for instance 
 \cite{bacciotti,khalil} and references therein). 

Concerning the stabilization of ensembles, we mention two papers using this approach:
in \cite{altafini2007}, the author aims at stabilizing an  ensemble of interacting spins along a reference trajectory; 
the result is achieved by showing, by means of 
Lie-algebraic methods, that the distance
between the state of the system and the target trajectory is a Lyapunov function. 
In \cite{ryan}, Jurdjevic-Quinn conditions  are applied to
stabilize an ensemble of harmonic oscillators. 

The feedback form of the control guarantees more robustness with respect to open-loop controls, and gives rise to a continuous bounded control, more easy to implement in
practical situations. 
We stress that, in the finite dimensional case,
the implementation of the control requires the knowledge of the \emph{bulk magnetization} of all spin, which is accessible through classical measurements 
(see for instance \cite{NMR-tanti,altafini2007}). 
We finally remark that the control proposed in this paper is very similar to the \emph{radiation damping effect} arising in NMR (see \cite{bloem,augustine}); we comment this fact
in the conclusion.

The structure of the paper is the following: in Section~\ref{sec:statement} we state the problem in general form; in Section~\ref{sec: finite} we tackle the finite dimensional case,
while in Section~\ref{sec: countable} we  analyze the case of a 
countable family of systems. Section~\ref{sec: simu} is devoted to some numerical results.

\section{Statement of the problem}\label{sec:statement}

We consider an ensemble of non-interacting spins immersed in a static magnetic field of strength $B_0(\boldsymbol{r})$, directed along the $z$-axis, and a time varying
transverse field $(B_x(t),B_y(t),0)$ (rf field), that we can control. The Bloch equation for this system takes then the form
\begin{equation} \label{eq: bloch}
\frac{\partial\boldsymbol{M}}{\partial t}(\boldsymbol{r},t)
=
\begin{pmatrix}
0 & -B_0(\boldsymbol{r}) & B_y(t)\\ 
B_0(\boldsymbol{r}) & 0 & -B_x(t)\\
-B_y(t) & B_x(t) & 0
\end{pmatrix}
\boldsymbol{M}((\boldsymbol{r}),t)
\end{equation}
(here for simplicity we set $\gamma=1$). For more details, we mention the monograph \cite{abragam}.

Since the dependence on the spatial coordinate $\boldsymbol{r}$ appears only in 
$B_0(\boldsymbol{r})$, we can represent $M(\boldsymbol{r},t)$ as a collection of time-dependent vectors $X_e(t)=(x_e(t) , y_e(t), z_e(t))$, 
where $e=B_0(\boldsymbol{r})$, each one  belonging to the unit sphere $S^2\subset \mathbb{R}^3$ and subject to the law

\begin{equation} \label{eq: system}
\begin{pmatrix}
\dot{x}_e \\ \dot{y}_e\\ \dot{z}_e
\end{pmatrix}
=
\begin{pmatrix}
0 & -e & u_2\\ 
e & 0 & u_1\\
-u_2 &  -u_1 & 0
\end{pmatrix}
\begin{pmatrix}
x_e \\ y_e\\ z_e 
\end{pmatrix},
\end{equation}
with $u_1(t)=-B_x(t)$ and $u_2(t)=B_y(t)$.
The Larmor frequencies $e$ of the spins in the ensemble take value in some subset $\mathcal{E}\subset I$ of a bounded  
interval $I$.
Depending on the spatial distribution of the spins, $\mathcal{E}$ could be a finite set, an infinite countable set, or an interval.

We are concerned with the following control problem:

\medskip
\noindent
{\bf (P)}{\it Design a control function $\uu:[0,+\infty) \to \mathbb{R}^2$ such that for every $e\in \mathcal{E}$ the solution of equation 
\eqref{eq: system} 
is driven to $X_e=(0, 0,-1)$.}
\medskip

To face this problem, 
we consider the Cartesian product $\Sf=\prod_{e\in \mathcal{E}} S^2$, whose elements are the collections 
$\bX=\{X_e\}_{e\in \mathcal{E}}$ such that $X_e\in S^2$ for every $e\in \mathcal{E}$. Depending on 
the structure of $\mathcal{E}$, $\bX$ can be a finite or an infinite countable collection of states $X_e\in S^2$, or a function $\bX : \mathcal{E} \to S^2$ belonging to some functional space.
The collection $\bX$ of magnetic moments evolves according to the equation
\begin{equation} \label{global Bloch}
 \dot{\bXi}=\bF(\bXi,\uu), \qquad \bXi(0)=\bX,
\end{equation}
where $\bF$ denotes the collection $\bF=\{F_{e}\}_{e\in \mathcal{E}}$ of tangent vectors to $S^2$, with $F_e(\bX,\uu)=\begin{pmatrix}
0 & -e & u_2\\ 
e & 0 & u_1\\
-u_2 &  -u_1 & 0
\end{pmatrix} X_e$, and $\uu=(u_1,u_2)$.

Some remarks on the existence of solutions for equation \eqref{global Bloch} are in order, and will be provided case by case. 
Assuming that these issues are already fixed, we define the two states $\bX^+=\{X_e : \forall e\in \mathcal{E} \ X_e=(0,0,1)\}$ and
$\bX^-=\{X_e :  \forall e\in \mathcal{E} \ X_e=(0,0,-1) \}$, and rewrite the problem {\bf (P)} as 

\medskip
\noindent
{\bf (P')}{\it  Design a control function $\uu:[0,+\infty) \to \mathbb{R}^2$ such that the solution of equation \eqref{global Bloch}  
is driven to $\bX=\bX^-$.}
\medskip

We remark that the notion of convergence of $\bX(\cdot)$ towards $\bX^-$ in problem {\bf (P')} has to be specified case by case, depending on the structure of the set $\mathcal{E}$ and on the
topology of $\Sf$.

\section{Finite dimensional case} \label{sec: finite}

First of all, we consider the case in which the set $\mathcal{E}$ is a finite collection of pairwise distinct energies, 
that is $\mathcal{E}=(e_1,\ldots,e_p)$ such that $e_k\in I \ \forall k\in\{1,\ldots,p\}$ and $e_k\neq e_j$ if $i\neq j$.
We recall that the state space $\Sf$ of the system is the finite product of $p$ copies of $S^2$.

\begin{lemma}
Assume that all energy levels $e_i$ are pairwise distinct. Let $\mathcal{I}=\{\bX\in \Sf: x_{e_i}=y_{e_i}=0\ \forall\, i=1,\ldots,p\}$.
Then every solution of the the control system \eqref{eq: system} with control  
\begin{equation} \label{eq: control}
\begin{cases}
u_1=\sum_{i=1}^p y_{e_i}
\\u_2=\sum_{i=1}^p x_{e_i}  
\end{cases}
\end{equation}
tends to $\mathcal{I}$ as $t\to +\infty$.
\end{lemma}

\proof
Consider the function $V(\bX)=\sum_{i=1}^p z_{e_i}$, and let $\bXi(\cdot)$ be a solution of \eqref{global Bloch} with the control given in \eqref{eq: control}. 
We notice that $\dot{V}(\bXi(t))=-\left(\sum_{i=1}^p x_{e_i}\right)^2-\left(\sum_{i=1}^p y_{e_i}\right)^2$, therefore it is non-positive on the whole $\Sf$,
and it is zero only on  the set $\mathcal{M}=\{\bX\in \Sf: \sum_{i=1}^px_{e_i}=\sum_{i=1}^py_{e_i}=0\}$.
We can then apply La Salle invariance principle to conclude that, for every initial condition, $\bXi(t)$ tends to the largest invariant subset of  $\mathcal{M}$.

Consider a trajectory $\bXi(\cdot)$ entirely contained in $\mathcal{M}$. Since $\uu=0$, then for every $i$ we have that
\[\bXi_i(t)=\begin{pmatrix}
\cos(e_i t) & -\sin(e_i t) & 0\\            
\sin(e_i t) & \cos(e_i t) & 0\\            
0 & 0 & 1
\end{pmatrix}
\begin{pmatrix}
x_{e_i}(0)\\ y_{e_i}(0) \\ z_{e_i}(0) 
\end{pmatrix}.
\]
By definition, for every $t\geq 0$ it holds $\sum_{i=1}^p x_{e_i}(t)=\sum_{i=1}^p y_{e_i}(t)=0$. Differentiating these equalities $p-1$ times and evaluating
at $t=0$ we obtain the two conditions
\[
\begin{pmatrix}
1 & 1& \ldots & 1\\ 
e_1 & e_2& \ldots & e_p\\
\vdots & \vdots & &\\
e_1^{p-1}& e_2^{p-1}& \ldots& e_p^{p-1}
\end{pmatrix}
\begin{pmatrix}
x_1(0)\\ x_2(0) \\ \vdots\\ x_p(0)
\end{pmatrix}
=
\begin{pmatrix}
0\\ 0 \\ \vdots\\ 0
\end{pmatrix}
\qquad
\begin{pmatrix}
1 & 1& \ldots & 1\\ 
e_1 & e_2& \ldots & e_p\\
\vdots & \vdots & &\\
e_1^{p-1}& e_2^{p-1}& \ldots& e_p^{p-1}
\end{pmatrix}
\begin{pmatrix}
y_1(0)\\ y_2(0) \\ \vdots\\ y_p(0)
\end{pmatrix}
=
\begin{pmatrix}
0\\ 0 \\ \vdots\\ 0
\end{pmatrix}.
\]
The determinant of the Vandermonde matrix here above is given by $\prod_{1\leq i<j\leq p}(e_i-e_j)$, which is non-zero under the assumptions. Therefore the two equations are satisfied if and only if 
$\bXi(0)\in \mathcal{I}$. It is immediate to see that $\mathcal{I}$ is the largest invariant subset of $\mathcal{M}$.  
\eproof

\medskip
The set $\mathcal{I}$ is composed by a collection of $2^p$ isolated points $\bQ^k=(Q^k_{1},\ldots,Q^k_p)$, $k=1,\ldots,2^p$, where $Q^k_j=(0,0,\alpha^k_j)$ and $|\alpha^k_j|=1$.
These points are equilibria for the  closed-loop system \eqref{eq: system}-\eqref{eq: control}.
We distinguish three cases:
\begin{itemize}
 \item if $\alpha^k_j=-1$ for every $j$, then $\bQ^k=\bX^-$, and  it is an asymptotically stable equilibrium for the system;
 \item if $\alpha^k_j=1$ for every $j$, then $\bQ^k=\bX^+$ is an unstable equilibrium for the system; in particular, $\bX^+$ is a repeller;
 \item all other points in $\mathcal{I}$ are neither attractor neither repellers, since each of these points is a saddle-point of $V$.
 \end{itemize}
 
These facts will be proved in the next section (see Proposition~\ref{hyperbolicity}, Lemma~\ref{real part 1} and Remark~\ref{real part 2}). 
We end this one recalling the following property of the basin of attraction; even though it is a standard result, we are providing a sketch of the proof.
 
\begin{lemma}
Let $\mathcal{B}$ be the basin of attraction of $\bX^-$. Then $\mathcal{B}$ is an open neighborhood of $\bX^-$ and, in the case $p>1$, there exists at least one $\bQ\in \mathcal{I}\setminus
\{\bX^+,\bX^-\}$ such that $\bQ\in \partial \mathcal{B}$.
\end{lemma}
\proof
Let us denote with $\phi^t$ the map that associates with each $\bX_0\in \Sf$ the solution at time $t$ of the control system \eqref{global Bloch}-\eqref{eq: control} with initial
condition equal to $\bX_0$.
By definition,
$\bX_0 \in \mathcal{B}$ if for every neighborhood $\mathcal{V}$ of $\bX^-$ there exists a time $T_{\mathcal{V},\bX_0}$ such that $\phi^t(\bX_0) 
\in \mathcal{V}$ for every $t> T_{\mathcal{V},\bX_0}$. By definition, $\mathcal{B}$ is $\phi^t$-invariant.

The asymptotic stability of $\bX^-$ and the continuous dependence of $\phi^t$ from initial conditions imply that 
 $\mathcal{B}$ is an open neighborhood of $\bX^-$. 

Consider now a point $\bX\in \partial \mathcal{B}$. Then there exists some $\bQ \in \mathcal{I}\setminus \{\bX^-, \bX^+\}$ such that $\phi^t(\bX)\to \bQ$ as $t\to +\infty$.
Let us fix $\epsilon>0$ and choose $\bar t>0$ such that $|\phi^t(\bX)-Q|<\epsilon/2 \ \forall t\geq \bar t$. 
By continuity with respect to initial conditions, there exists $\delta>0$ such that if $|\bX-\bX'|<\delta$, then $|\phi^{\bar t}(\bX)-\phi^{\bar t}(\bX')|< \epsilon/2$, which 
implies that
 $|\phi^t(\bX')-\bQ|<\epsilon$, that is $\bQ\in \partial \mathcal{B}$. 
\eproof

\subsection{Linearized system}

In order to study the structure of the basin of attraction $\mathcal{B}$, we linearize the system \eqref{global Bloch}-\eqref{eq: control} around a 
point $\bQ\in\mathcal{I}$, and we study the corresponding eigenvalues. 
We will show below that the linearized system is always hyperbolic (when we consider its restriction to the tangent space to the collection of spheres).

The linearization gives
\begin{equation} \label{eq: linearised}
\begin{pmatrix}
\dot{\boldsymbol{\delta x}}\\\dot{\boldsymbol{\delta y}}  \\\dot{\boldsymbol{\delta z}} 
\end{pmatrix}=
\begin{pmatrix}
K_{\bQ} & -E & 0\\
E & K_{\bQ} & 0\\
0 & 0 & 0
\end{pmatrix}
\begin{pmatrix}
\boldsymbol{\delta x}\\\boldsymbol{\delta y} \\\boldsymbol{\delta z} 
\end{pmatrix}
\end{equation}
where 
\[
K_{\bQ}=
\begin{pmatrix}
z_1^{\bQ} & z_1^{\bQ} & \ldots & z_1^{\bQ}\\ 
z_2^{\bQ} & z_2^{\bQ} & \ldots & z_2^{\bQ}\\ 
\vdots & \vdots &  & \\ 
z_p^{\bQ} & z_p^{\bQ} & \ldots & z_p^{\bQ}\\ 
\end{pmatrix}
\qquad
E=
\begin{pmatrix}
e_1 & 0 & \ldots & 0\\
0 & e_2 & \ldots & 0\\
\vdots & \vdots &  & \\ 
0 & 0 & \ldots & e_p\\ 
\end{pmatrix}
\]
and $z_i^{\bQ}$ is the value of the coordinate $z_i$ at the point $\boldsymbol{Q}$.
Set moreover $M_{\bQ}=\left(\begin{smallmatrix}
K_{\bQ} & -E\\
E & K_{\bQ}
\end{smallmatrix}\right)$.
Notice that we can write $K_{\bQ}=\kappa_{\bQ} \zeta^T$, where $\kappa_{\bQ}=(z_1^{\bQ},\ldots,z_p^{\bQ})$ and $\zeta=(1, 1, \ldots,1)$, then  
rank$K_{\bQ}=1$.

In the following, 
with a little abuse of notation, we will remove the dependence on $\bQ$ from $K$, $M$, $\kappa$ and its components, specifying it only when necessary.

In order to compute the eigenvalues of the matrix $M$, we consider the complexification of system \eqref{eq: linearised}, that is we set $\xi=\boldsymbol{\delta x}+i\boldsymbol{\delta y}$, observing that 
$\dot{\xi}=(K+iE)\xi$.
  It is easy to see that $\ell$ is an eigenvalue of $M$ if and only if it is also either an eigenvalue of $(K+iE)$ or an eigenvalue of $(K-iE)$, 
that is, the spectrum of $M$ is equal to the union of the spectra of $(K+iE)$ and $(K-iE)$.

\medskip
\noindent
\textbf{Properties} In the following, we will use the following properties of block matrices
\begin{description}
 \item[(P1)] Let $M$ be the block matrix $\left(\begin{smallmatrix}
A & B\\
C & D
\end{smallmatrix}\right)$. If $A$ is invertible, then 
$\det M =\det(A) \det(D-CA^{-1}B)$. If $D$ is invertible, then $\det M=\det (D) \det (A-BD^{-1}C)$.
 \item[(P2)] Let $A$ be an invertible matrix of size $n$, and $x,y$ two $n$-dimensional vectors. Then $\det(A+xy^T)=\det(A) (1+y^TA^{-1}x)$.
\end{description}

\begin{lemma}
The matrices $(K+iE)$ and $(K-iE)$ are invertible.
\end{lemma}
\proof
Assume that $E$ is invertible.
Then $\det(K+iE)=\det(iE)\det(\mathbb{1}-iE^{-1}K)$, where $\mathbb{1}$ denotes the $n$-dimensional identity matrix; since $-iE^{-1}K=-iE^{-1}\kappa\zeta^T=\widetilde{\kappa}\zeta^T$, 
by {\bf (P2)} we have that 
$\det(\mathbb{1}+\widetilde{\kappa}\zeta^T)=(1+\zeta^T\widetilde{\kappa})\neq 0$, since $\zeta^T\widetilde{\kappa}$ is purely imaginary. 

If $E$ is not invertible, up to permutations and relabeling we assume that $e_1=0$. We suitably add or subtract the first row of $K$ to all other ones, in order to get that
\[
\det (K+iE)=\det
\begin{pmatrix}
z_1 & z_1 & \ldots &  z_1\\
0 & ie_2 & \ldots & 0\\
\vdots & & \ddots & \vdots\\
0 & 0 &\ldots & ie_p\\
\end{pmatrix}=(i)^{p-1}z_1 e_2\cdots e_p.
\]

The same arguments prove that 
$(K-iE)$ is invertible.
\eproof


\begin{proposition} \label{hyperbolicity}
For every $\bQ\in \mathcal{I}$, all eigenvalues of the matrix $M_{\bQ}$ have non-zero real part.
\end{proposition}

\proof
First of all, we prove by contradiction that  $ie_i$ is not an eigenvalue of $K+iE$ and $-ie_i$ is not an eigenvalue of $K-iE$, for every $i=1,\ldots,p$.
Assume that $i e_i$ is an eigenvalue of $K+iE$, that is there exists a vector $v\in \mathbb{C}^p$ such that $(K+iE)v=ie_iv$, that is 
$Kv=(i(e_i-e_1)v_1,\ldots,i(e_i-e_p)v_p)^T$; since $Kv=\kappa(\zeta^Tv)$ and all the components of $\kappa$ are different from zero,
this implies that $\zeta^Tv=0$ and therefore $Kv=0$. Then $i(E-e_i\mathbb{1})v=0$, that is $v_j=0$ for every $j\neq i$, therefore
$(Kv)_l=z_l\sum_{k=1}^p v_k=z_l v_i$ for every $l$. Since $Kv=0$, then $v=0$. Then $ie_i$ cannot be an eigenvalue of $K+iE$.
The corresponding statement for $K-iE$ is proved using an analogous argument.

Let us now assume, by contradiction, that $i\ell$, $\ell\in \mathbb{R}$, is an eigenvalue of $M$ relative to the eigenvector $(X,Y)$ (where $X,Y\in \mathbb{C}^p$), and assume that $X+iY\neq 0$ 
(if this is not the case, then  $X-iY\neq 0$ and we can repeat the same argument used below
with $(K-iE)$). Then $X+iY$ is an eigenvector of $(K+iE)$ relative to $i\ell$, and $\ell$ is different from every of the $e_i$. 
Let  $\chi,\eta \in \mathbb{R}^p$ be the real vectors such that $X+iY=\chi+i \eta$.
Straight computations show that
\[
\begin{cases}
K\chi-E\eta=-\ell \eta\\
K\eta+E\chi=\ell \chi
\end{cases}\Rightarrow
\begin{cases}
(E-\ell \mathbb{1}) \eta=K\chi=(\zeta^T \chi)\kappa\\
(E-\ell \mathbb{1}) \chi=-K\eta=-(\zeta^T \eta)\kappa.
\end{cases}
\]
Since $(E-\ell \mathbb{1})$ is invertible, we have that 
\begin{align*}
\chi&=- (\zeta^T \eta) (E-\ell \mathbb{1})^{-1}\kappa \\
\eta&= (\zeta^T \chi) (E-\ell \mathbb{1})^{-1}\kappa,
\end{align*}
that is, $\eta$ and $\chi$ are parallel and $X+iY=(a+ib)(E-\ell \mathbb{1})^{-1}\kappa$, for some real coefficients $a,b$. 
Then $(E-\ell \mathbb{1})^{-1}\kappa$ is a real eigenvector of 
$(K+iE)$ relative to $i\ell$, which implies that $K(E-\ell \mathbb{1})^{-1}\kappa=0$ and $E(E-\ell \mathbb{1})^{-1}\kappa=\ell(E-\ell \mathbb{1})^{-1}\kappa$, which is possible
only if 
$\kappa$ is null.
\eproof

\begin{lemma} \label{real part 1}
Let $\ell$ be an eigenvalue of $K_{\bQ}+iE$. Then the following equality holds
\begin{equation} \label{eigenvalue fin}
\sum_{j=1}^p \frac{z_j^{\bQ} (\lambda+i (\mu-e_j))}{\lambda^2+(e_j-\mu)^2}=1, 
\end{equation}
where $\lambda$ and $\mu$ denote respectively the real and the imaginary part of $\ell$.
In particular, all the eigenvalues of $K_{\bX^+}+iE$ have positive real part and all 
the eigenvalues of $K_{\bX^-}+iE$ have negative real part.
\end{lemma}

\proof
Thank to property {\bf (P2)} and the fact that $\ell$ is not an eigenvalue of $iE$, it holds
\begin{align*}
\det (K+iE-\ell \mathbb{1}) &= \det(iE-\ell \mathbb{1})(1+\zeta^T (iE-\ell \mathbb{1})^{-1}\kappa)\\
&=\det(iE-\ell \mathbb{1})\Big(1+\sum_{j=1}^p \frac{z_j}{ie_j-\ell}\Big).
\end{align*}
Equation \eqref{eigenvalue fin} follows from $1+\sum_{j=1}^p \frac{z_j}{ie_j-\ell}=0$.

In particular, since $z^{\bX^-}_i=-1$ for every $i$, for every eigenvalue $\ell$ of $K_{\bX^-}+iE$ equation \eqref{eigenvalue fin} reads  
\[\begin{cases}
\sum_{j=1}^p \frac{ \lambda}{\lambda^2+(e_j-\mu)^2}=-1 \\ 
\sum_{j=1}^p \frac{ (\mu-e_j)}{\lambda^2+(e_j-\mu)^2}=0
\end{cases},
\]
which implies $\lambda<0$. The same argument proves that all eigenvalues of $K_{\bX^+}+iE$ have positive real part.
\eproof

Analogous computations show that every eigenvalue $\ell=\lambda+i\mu$ of $K_{\bQ}-iE$ satisfies  the equation 
\[\sum_{j=1}^p \frac{z_j^{\bQ} (\lambda-i (\mu-e_j))}{\lambda^2+(e_j-\mu)^2}=1,\]
and that all  eigenvalues of $K_{\bX^-}-iE$ (respectively, $K_{\bX^+}-iE$)  have negative (respectively, positive) real part.

\medskip

Let us now consider the linearized flow in the tangent space to $\Sf$ at some $\bQ\in \mathcal{I}$. First of all, we notice that $T_{\bQ}\Sf=\{(\boldsymbol{\delta x},
\boldsymbol{\delta y}, 0)\}$ for every $\bQ\in \mathcal{I}$, therefore the linearization of the flow $\phi^t$ on $T_{\bQ} \Sf$ can be represented by the matrix $M$. In particular, 
Proposition~\ref{hyperbolicity} implies that each $\bQ\in \mathcal{I}$ is a hyperbolic equilibrium for the flow $\phi^t$ (restricted to $\Sf$).

\begin{remark} \label{real part 2}
For every $\bQ\neq \bX^-$, then at least two eigenvalues of $D_{\bQ}\phi^t$ have positive real part. Indeed, Proposition~\ref{hyperbolicity} implies that all 
eigenvalues of $D_{\bQ}\phi^t$
have non-zero real part; if they all had negative real part, there would be a 
contradiction with the fact none of the $\bQ \in \mathcal{I}\setminus\{\bX_-\}$ is a local minimum of the Lyapunov function, that is, none of these equilibria is stable.
Since the eigenvalues of $D_{\bQ}\phi^t$ come in conjugate pairs\footnote{for purely real eigenvalues, this means that they must have even multiplicity}, at least two must have positive real part.
\end{remark}

We are now ready to state of the main results.

\begin{theorem} \label{th: 1}
The point $\bX^-$ is asymptotically stable for the system \eqref{eq: system}-\eqref{eq: control}, and its basin of attraction is open and dense in $\Sf$.
\end{theorem}

\proof
For $p=1$, then the basin of attraction of $\bX^-$ is trivially $S^2\setminus \bX^+$. Let us then assume that $p>1$.

Consider $\bQ \in \mathcal{I}\setminus\{\bX^+,\bX^-\}$. From Proposition \ref{hyperbolicity}, Lemma~\ref{real part 1} and Remark~\ref{real part 2}, we know that the restriction of $D_{\bQ}\phi^t$ to $T_{\bQ}\Sf$
 satisfies the following properties: there exists a splitting of the tangent space $T_{\bQ}\Sf=E_{\bQ}^-\oplus E_{\bQ}^+$ such that
 \begin{itemize}
  \item there exists $\rho_+>1$ such that $\|D\phi^{-t} |_{E_{\bQ}^+}\|\leq \rho_+^{-t}$ and $\dim E_{\bQ}^+\geq 2$
 \item there exists $\rho_-<1$ such that $\|D\phi^{t} |_{E_{\bQ}^-}\|\leq \rho_-^t$ and $\dim E_{\bQ}^-\geq 2$. 
 \end{itemize}
Then we can apply Hadamard-Perron Theorem \cite{Katok-Hasselblatt} and conclude that there exist two $\mathcal{C}^1$-smooth
injectively immersed submanifolds 
$W_{\bQ}^s,W_{\bQ}^u\subset \Sf$ such that 
\begin{gather}
W_{\bQ}^s=\{X\in \Sf : \mathrm{dist}(\phi^t(X), \bQ)\to 0 \ \mathrm{as} \ t\to +\infty\} \quad \mathrm{and} \quad T_{\bQ}W_{\bQ}^s=E_{\bQ}^-\\
W_{\bQ}^u=\{X\in \Sf : \mathrm{dist}(\phi^{-t}(X), \bQ)\to 0 \ \mathrm{as} \ t\to +\infty\} \quad \mathrm{and} \quad T_{\bQ}W_{\bQ}^u=E_{\bQ}^+.
\end{gather}
We recall that every point in $\Sf$ asymptotically reaches $\mathcal{I}$, under the action of the flow $\phi^t$. Therefore,
the set of all points that do not asymptotically reach $\bX^-$ is
\[\mathcal{Q}=\bigcup_{\bQ\in \mathcal{I}\setminus \{\bX^+, \bX^- \}} W_{\bQ}^s \cup \{\bX^+\}.\] 

Set $\bG^p=\Sf \setminus \mathcal{Q}$, 
and notice that $\mathcal{Q}$ is a finite union of smooth manifolds of codimension at least $2$. 
This implies that its complement is dense.

Let $\bX^k$ be a sequence in $\mathcal{Q}$, converging to some $\bar{\bX}\in \Sf$. By continuity with respect to initial conditions, for every $\epsilon>0$ and every 
$T>0$ there exists $\bar k$
such that if $k\geq \bar{k}$, then $|\phi^T(\bX^k)-\phi^T(\bar{\bX})|\leq \epsilon$, which implies, by smoothness 
of the Lyapunov function $V$, 
that $|V(\phi^T(\bX^k))-V(\phi^T(\bar{\bX}))|\leq L\epsilon$, for some $L>0$. 

Since for every $t$ and every $k$ it holds 
$V(\phi^t(\bX^k))\geq \bar{V}$, where $\bar{V}=\min_{\mathcal{I}\setminus \{\bX^-\} } V$,   then we can conclude that for every $\epsilon>0$ and  $T>0$ we can find $\bar k$
such that 
$V(\phi^T(\bar{\bX}))\geq V(\phi^T(\bX_{\bar k}))-\epsilon\geq \bar{V}-\epsilon$. Then 
\[\phi^t(\bar{\bX})\to  \mathcal{I}\setminus \{\bX^-\}  ,\]
that is
\[
\overline{\mathcal{Q}}\subset \bigcup_{\bQ\in \mathcal{I}\setminus \{\bX^-\}}W_{\bQ}^s=\mathcal{Q}.
\]
\eproof

\begin{remark}
Theorem~\ref{th: 1} states that the set of ``bad'' initial conditions - that is, the set of initial condition not converging to the state $\bX^-$ - is given by the union of the unstable
equilibria and of their corresponding stable manifold. In the single spin case, the ``bad set'' reduces to the stable equilibrium $\bX^+$, as already pointed out in 
\cite{altafini2007}, where a similar feedback control is applied for stabilizing a set of interacting spins.
\end{remark}

\section{Countable case}\label{sec: countable}

\subsection{Existence of solutions}

Let us now assume that $\mathcal{E}=\{e_i\}_{i\in \mathbb{N}^+}$ is a sequence of pairwise distinct elements contained in $I$.
The state of the system is represented by the sequence $\bX=\{X_{e_i}\}_i$, with $X_{e_i} \in S^2$, and  
the state space is the countable Cartesian product $\Sf=\Pi_{i=1}^{\infty}S^2$.

Before trying to solve the problem {\bf (P')}, it is necessary to discuss its well-posedness. To do this,
let us consider the function $\bd$ on the infinite Cartesian product $\Pi_{i=1}^{\infty}\mathbb{R}^3$: 
\[
\bd(\bX,\bX')=\sum_{i=1}^{\infty} w_i |X_{e_i}-X_{e_i}'|,  
\]
where $|X_{e_i}-X_{e_i}'|^2=|x_{e_i}-x_{e_i}'|^2+|y_{e_i}-y_{e_i}'|^2+|z_{e_i}-z_{e_i}'|^2$ and $\{w_i\}_{i\in \mathbb{N}^+}$ is a positive monotone 
sequence 
such that the series $\sum_{i\in \mathbb{N}^+}w_i$ converges. Without loss of generality, here and below we put $w_i=2^{-i}$.
We now consider the subset $\Pib\subset \Pi_{i=1}^{\infty}\mathbb{R}^3$ of all sequences $\bX$ such that $\sum_{i=1}^{\infty} w_i |X_{e_i}|$ is finite. 

It is immediate to see that $\Pib$, endowed with the distance function $\bd$, is a Banach space.
More precisely, it corresponds to a weighted $\ell_1$-space. The choice of a weighted space is motivated by two exigencies: first of all, 
it guarantees the compactness of the set  $\Sf$ with respect to the topology induced by $\bd$, as will be proved in the next section; in addition, it permits to define the Lyapunov function and
the feedback control as straightforward extensions of those of the finite-dimensional case, avoiding well-definiteness issues.

We remark that $\Sf$ is a proper connected subset of the unit sphere in the Banach space $(\Pib,\bd)$.

\begin{remark}
By standard arguments, it is easy to prove that, for every  $-\infty<a<b<+\infty$,   $\mathcal{C}([a,b],\Pib)$ is a Banach space with respect to the sup norm
\[
\|\boldsymbol{f}\|_{\mathcal{C}([a,b],\Pib)} =\sup_{t\in[a,b]} \Big|\sum_{i=1}^{\infty} 2^{-i} |f_i(t)|\Big|.
\]
\end{remark}

We now consider the feedback control $\uu=(u_1,u_2)$, defined by
\begin{equation} \label{eq: control 2}
\begin{cases}
u_1=\sum_{i=1}^{\infty} 2^{-i} y_{e_i}
\\u_2=\sum_{i=1}^{\infty}2^{-i} x_{e_i}  
\end{cases}
\end{equation}
and we plug it into the control system \eqref{global Bloch}. 
The resulting autonomous dynamical system on $\Sf$ is well defined, as the following result states.

\begin{theorem}
The Cauchy problem $\dot{\bXi}=\boldsymbol{F}(\bXi,\uu)$ with initial condition in $\Sf$ is well-defined. 
\end{theorem}

\proof
In order to apply the standard existence theorem of solution of ODEs in Banach spaces (see for instance \cite{amann,deimling}), we need our solution space 
to be a linear space. Therefore, we consider the Cauchy problem on $\Pib$
\begin{equation}
\begin{cases} \label{eq: system cut-off} 
\dot{\bXi}=\widetilde{\boldsymbol{F}}(\bXi)\\
\bXi(0)=\bXi^0,
\end{cases}
\end{equation}
where $\widetilde{\boldsymbol{F}}=\{\widetilde{F}_i\}_{i=1}^{\infty}$ is the vector field on $\Pib$ defined by
\[
\widetilde{F}_i(\bX)= e_i A \psi(X_{e_i})+\varphi\Big(\sum_{j=1}^{\infty}2^{-j} x_{e_j}\Big)B\psi(X_{e_i})+\varphi\Big(\sum_{j=1}^{\infty}2^{-j} y_{e_j}\Big)C\psi(X_{e_i}), 
\]
where
$
A=\left(
\begin{smallmatrix}
0 & -1 & 0\\
1 & 0 & 0\\
0 & 0 & 0
\end{smallmatrix}\right)
$,
$B=\left(
\begin{smallmatrix}
0 & 0 & 1\\
0 & 0 & 0\\
-1 & 0 & 0
\end{smallmatrix}\right)$,
$C=\left(
\begin{smallmatrix}
0 & 0 & 0\\
0 & 0 & -1\\
0 & 1 & 0
\end{smallmatrix}\right)$,
and $\varphi :\mathbb{R}\to\mathbb{R}$ and $\psi :\mathbb{R}^3\to\mathbb{R}^3$ are the cut-off functions
\[
\varphi(x)=
\begin{cases}
x & \mathrm{if} \ |x|\leq b\\
b & \mathrm{if} \ |x|\geq b
\end{cases}
\qquad
\psi(\ww)=
\begin{cases}
\ww & \mathrm{if}\ |\ww|\leq a\\
a\frac{\ww}{|\ww|} & \mathrm{if}\ |\ww|\geq a
\end{cases},
\]
for some real numbers $a,b>1$.  The uniform Lipschitz continuity of $\widetilde{\boldsymbol{F}}$ can be easily proved by computations.
Applying the Picard-Lindel\"of Theorem (\cite{amann,deimling}), we obtain
that there exists an interval $I_0$ containing 0 such that
the Cauchy problem \eqref{eq: system cut-off} admits a unique solution, continuously differentiable on $I_0$.

We notice that, if the initial condition belongs to $\Sf$, then the solution of \eqref{eq: system cut-off} belongs to $\Sf$ for all $t\in I_0$. Moreover, 
by the global Lipschitz continuity of $\widetilde{\boldsymbol{F}}$, we deduce that  
 the solution arising from any initial condition in $\Sf$  is well defined for all $t\in \mathbb{R}$.
Finally, we observe that
$\widetilde{\boldsymbol{F}}|_{\Sf}=\boldsymbol{F}$, therefore the solutions of \eqref{eq: system cut-off} with initial condition in $\Sf$ coincide with the solutions of
the equation $\dot{\bXi}=\bF(\bXi,\mathbf{u}(\bXi))$ with the same initial condition. 
\eproof

A direct application of Gronwall inequality yields the following result.

\begin{proposition}
The solutions of the Cauchy problem \eqref{eq: system cut-off} depend continuously on initial conditions. 
\end{proposition}

\subsection{Asymptotic pointwise convergence to $\bX^-$}

Let us consider the function $V: \Sf \to \mathbb{R}$ defined by $V(\bX)=\sum_{i=1}^{\infty} 2^{-i}z_{e_i}$. It is easy to see that its time derivative along the integral curves 
of the vector field $\bF(\bX,\mathbf{u}(\bX))$ satisfies $\dot{V}=-\big( \sum_{i=1}2^{-i} x_{e_i} \big)^2-\big( \sum_{i=1}2^{-i} y_{e_i} \big)^2 \leq 0$. In order to conclude about the stability 
of these trajectories by means of a La Salle-type argument, we need to prove that $\Sf$ is compact.
To do that, let us first recall the following definition (see for instance \cite{kelley}).
\begin{definition}
The product topology $\mathcal{T}$ on $\Sf$ is the coarsest topology that makes continuous all the projections $\pi_i : \Sf \to S_{e_i}^2$.
\end{definition}

By Tychonoff's Theorem, any product of compact topological spaces is compact with respect to the product topology  (\cite{kelley}). This in particular implies that
$\Sf$ is compact with respect to $\mathcal{T}$.

As we will see just below, the product topology is equivalent to topology induced by the distance $\bd$, so $\Sf$ is compact with respect to the latter.
\begin{lemma}
Let us denote with $\mathcal{T}_{\bd}$ the topology on $\Sf$ induced by $\bd$. We have that $\mathcal{T}=\mathcal{T}_{\bd}$.
\end{lemma}

\proof
By definition, $\mathcal{T}\subset\mathcal{T}_{\bd}$. If we prove that the open balls (that are a basis for $\mathcal{T}_{\bd}$) are open with respect to $\mathcal{T}$,
then $\mathcal{T}_{\bd}\subset\mathcal{T}$ and we get the result.

Let $N>0$ and let us define the function $\bd^N: \Sf \times \Sf \to [0,1]$  as
$\bd^N(\bX,\bX')=\sum_{i=1}^N 2^{-i} |\bX_i-\bX_i'| $. It is easy to prove that $\bd^N$ is continuous with respect to $\mathcal{T}$; indeed, the restriction  
$\bd|_{\prod_{i=1}^N S_i^2}$ is obviously continuous with respect to the product topology on $\prod_{i=1}^N S_{e_i}^2$,
and for every open interval $(a,b)\subset[0,1]$ we have that 
$(\bd^{N})^{-1}(a,b)=\big(\bd|_{\prod_{i=1}^N S_{e_i}^2}\big)^{-1}(a,b)\times \prod_{i>N}S_{e_i}^2$,
which is open with respect to $\mathcal{T}$.

The sequence $\{\bd^N\}^N$ converges uniformly to $\bd$. Indeed, for every $\bX,\bX'\in \Sf$ we have that
\[
|\bd^N(\bX,\bX')-\bd(\bX,\bX')|=|\sum_{k\geq N+1} 2^{-i} |X_{e_i}-X_{e_i}'||\leq 2^{-N} .
\]
Then $\bd$ is continuous with respect to $\mathcal{T}$, and this completes the proof.
\eproof

Thanks to previous Lemma, we can conclude that $\Sf$ is compact with respect to $\bd$. 
In particular, this permits to prove a version of La Salle invariance principle holding for 
the equation $\dot{\bXi}=\bF(\bXi,\uu(\bXi))$ (the result can be found for instance in \cite[Theorem 18.3, Corollary 18.4]{amann}; nevertheless, for completeness
in the exposition, we are giving a proof here below).

We also remark that the Lyapunov function $V$ is continuous with respect to $\mathcal{T}$.

\begin{proposition}[Adapted La Salle] \label{adapted LaSalle}
Let us consider the set $\mathcal{M}=\{\bX \in \Sf : \dot{V}(\bX)=0\}$, where we use the notation $\dot{V}(\bX)=\frac{d}{dt} V(\phi^t(\bX)|_{t=0}$, and
$\phi^t$ denotes the flow associated with the dynamical system $\dot{\bXi}=\bf(\bXi,\uu(\bXi))$. Let $\mathcal{I}$
be the largest subset of $\mathcal{M}$ which is invariant for the flow $\phi^t$. Then for every $\bX\in \Sf$ we have that
$\phi^t(\bX)\to \mathcal{I}$ as $t\to+\infty$.
\end{proposition}

\proof 
The proof of this proposition relies on the compactness of $\Sf$ with respect to the topology $\mathcal{T}_{\boldsymbol{d}}$, and follows standard arguments.

Let $m=\min_{\bX\in\Sf}V(\bX)$, and 
fix $\bX^0\in \Sf$. By continuity of $V$, there exists $a=\lim_{t\to+\infty}V(\phi^t(\bX^0))$, $a\geq m$.

Let $\Omega_{\bX^0}=\{\bX\in \Sf : \exists \ (t_n)_n\to +\infty : \phi^{t_n}(\bX^0)\to\bX\}$ denote the $\omega$-limit set issued from $\bX_0$; notice that $\Omega_{\bX^0}$ is non-empty, since $\Sf$ is compact, therefore
for every sequence $(t_n)_n \to +\infty$ there exists a subsequence $(t_{n_k})_k$ such that $\phi^{t_{n_k}}(\bX^0)$ converges to some point in $\Sf$.
It is easy to see that $\Omega_{\bX^0}$ is invariant for the flow $\phi^t$ and therefore, since by continuity  $V|_{\Omega_{\bX^0}}=a$, we obtain  that $\dot{V}(\bX)=0$ 
for every $\bX\in \Omega_{\bX^0}$. This implies that $\Omega_{\bX^0}\subset\mathcal{M}$.

Let us now prove that $\Omega_{\bX^0}$ is compact. Consider a sequence $(\bX^k)_k$ contained in $\Omega_{\bX^0}$; by compactness of  $\Sf$, it converges, up to subsequences, 
to some $\bar{\bX}\in\Sf$ (we relabel the indexes).
By definition, for every $k$ there exist a sequence $(t_{k_n})_n\to+ \infty$ such that $\lim_{n}\phi^{t_{k_n}}(\bX^0)=\bX^k$. Moreover,
it is possible to define a divergent sequence $(\tau_k)_k$ such that
$\boldsymbol{d}(\phi^{\tau_k}(\bX_0),\bX^k)\leq 1/2k$ for every $k$. 
Fix $\epsilon>0$ and choose some $\bar k\geq1/\epsilon$ such that $\boldsymbol{d}(\bar{\bX},\bX^k)\leq \epsilon/2$ for $k\geq \bar{k}$ (possibly taking a suitable subsequence).
Then for $k\geq \bar{k}$ we have that $\boldsymbol{d}(\bar{\bX},\phi^{\tau_k}(\bX^0))\leq 
\epsilon$.
This means that $\bar{\bX}\in \Omega_{\bX^0}$, that is $\Omega_{\bX^0}$ is compact. 

Finally, let us assume, by contradiction, that there exist an open neighborhood $U$ of $\Omega_{\bX^0}$ in $\Sf$ and a sequence $(t_n)_n\to +\infty$ such that 
$\phi^{t_n}(\bX^0)\in \Sf\setminus U$ for every $n$. By compactness of $\Sf$, $\phi^{t_n}$ converges up to subsequences to some $\bar{\bX} \in \Sf \setminus U$.
 But by definition $\bar{\bX}\in \Omega_{\bX^0}$, then we have a contradiction.
 
 Let us now set $\mathcal{I}=\cup_{\bX \in \Sf} \Omega_{\bX}$. By construction, it is an invariant subset contained in $\mathcal{M}$.
 \eproof

By definition, $\mathcal{M}=\{\bX \in \Sf : \sum_{i=1}^{\infty}2^{-i}x_{e_i}=\sum_{i=1}^{\infty}2^{-i}y_{e_i}=0 \}$. Now we look for its largest invariant subset.
Let $\bX^0\in \mathcal{M}$; with the same argument than above, we can see that $\phi^t(\bX^0)=\{\ (x_{e_i}(t),y_{e_i}(t),z_{e_i}(t))\}_i$ with
\[
x_{e_i}(t)=\cos(e_i t)x_{e_i}^0-\sin(e_i t)y_{e_i}^0\qquad y_{e_i}(t)=\sin(e_i t)x_{e_i}^0+\cos(e_i t)y_{e_i}^0 \qquad z_{e_i}(t)=z_{e_i}^0.
\]

If $\bX^0$ belongs to an invariant subset of $\mathcal{M}$, then $\sum_{i=1}^{\infty}2^{-i}x_{e_i}(t)=\sum_{i=1}^{\infty}2^{-i}y_{e_i}(t)=0$ for every $t$.
Let us consider the two functions 
\begin{align*}
f(t)&=\sum_{i=1}^{\infty}2^{-i}x_{e_i}(t)=\sum_{i=1}^{\infty} 2^{-i}(\cos(e_i t)x_{e_i}^0-\sin(e_i t)y_{e_i}^0)\\
g(t)&=\sum_{i=1}^{\infty}2^{-i}y_{e_i}(t)=\sum_{i=1}^{\infty} 2^{-i}(\sin(e_i t)x_{e_i}^0+\cos(e_i t)y_{e_i}^0).
\end{align*}
It is easy to see that both $f(t)$ and $g(t)$ are uniform limits of trigonometric 
polynomials, that is, they are almost periodic functions (also referred to as \emph{uniform almost periodic functions} or \emph{Bohr almost 
periodic functions}, see \cite{besicovitch,corduneanu}). According to the references \cite{besicovitch,corduneanu}, the Fourier series of a (uniform) almost periodic function 
is computed as follows: for every $\omega \in \mathbb{R}$, we define the Fourier coefficients as
\[
a(f,\omega)=\lim_{T\to \infty}\frac{1}{T}\int_0^T f(t) e^{-i\omega t}\;dt \qquad
a(g,\omega)=\lim_{T\to \infty}\frac{1}{T}\int_0^T g(t) e^{-i\omega t}\;dt.
\]
By easy computations, we see that $a(f,\omega)$ and $a(g,\omega)$ are zero for every $\omega\notin\{e_i,-e_i\}_{i\in \mathbb{N}^+}$, and, moreover, that
\[
a(f,e_i)=\frac{x_{e_i}^0+iy_{e_i}^0}{2^{i+1}}, \quad 
a(f,-e_i)=\frac{x_{e_i}^0-iy_{e_i}^0}{2^{i+1}}, \quad
a(g,e_i)=\frac{y_{e_i}^0-ix_{e_i}^0}{2^{i+1}}, \quad 
a(g,-e_i)=\frac{ix_{e_i}^0+y_{e_i}^0}{2^{i+1}}.
\]
The Fourier series of $f$ and $g$ are respectively 
\[
f(t) \sim \sum_{i=1}^{\infty}a(f,e_i)e^{ie_i t}+a(f,-e_i)e^{-ie_i t}
\quad
g(t) \sim \sum_{i=1}^{\infty}a(g,e_i)e^{ie_i t}+a(g,-e_i)e^{-ie_i t}.
\]

By 
\cite[Theorem 1.19]{corduneanu}, the functions $f$ and $g$ are identically zero if and only if all the coefficients in their Fourier series are all null,
that is $x_{e_i}^0=y_{e_i}^0=0$ for every $k$.
Then the largest invariant subset of $\mathcal{M}$ is $\mathcal{I}=\{\bX : x_{e_i}=y_{e_i}=0 , |z_{e_i}|=1 \ \forall \ i\}$.

Applying Proposition \ref{adapted LaSalle} to these facts, we get the following result.
\begin{corollary}
Let $\bX^0\in \Sf$, and let $\bXi(t)=\phi^t(\bX^0)$, with the usual notation $\bXi(t)=\{\Xi_{e_i}(t)\}_i$ and $\Xi_{e_i}(t)=(x_{e_i}(t),y_{e_i}(t),z_{e_i}(t))$.
Then 
\begin{gather*}
\lim_{t\to+\infty} x_{e_i}(t)= \lim_{t\to+\infty} y_{e_i}(t)=0\\
|\lim_{t\to+\infty} z_{e_i}(t)|=1
\end{gather*}
 for every $i\geq 1$.
\end{corollary}

\begin{remark}
It is easy to see that every point $\bX \in \mathcal{I}$ is an accumulation point for the set $\mathcal{I}$,
but $\mathcal{I}$ is not dense in $\Sf$.
In particular, $\bX^-$ is an accumulation point for $\mathcal{I}$ too.
\end{remark}

\medskip

In the following, for every $N\geq 1$ we consider the \emph{truncated feedback control} $\uu^N=(u_1^N,u_2^N)$, 
where $u_1^N=\sum_{i=1}^N 2^{-i} y_{e_i}$ and $u_2^N=\sum_{i=1}^N 2^{-i} x_{e_i}$, and
we call $\bXi^N(\cdot)$ the solution of the differential equation  $\dot{\bXi}^N=\boldsymbol{F}(\bXi^N,\uu^N)$,  $\bXi^N \in \mathfrak{X}$. We remark that, for initial conditions in $\Sf$, the solution of the
Cauchy problem exists  and remains in $\Sf$ for all $t$.
As above, we use the notations
$\bXi^N(\cdot)=\{\Xi^N_{e_i}(\cdot)\}_i$ with $\Xi^N_{e_i}=\{(x^N_{e_i}(\cdot),y^N_{e_i}(\cdot),z^N_{e_i}(\cdot))\}_i$.

Applying the same arguments as in Section \ref{sec: finite}, we can prove the following result.
\begin{proposition} \label{prop: trunc}
For every $N>0$, there exists an open dense set 
$A\subset \Sf$ such that for every $\bX\in A$ the solution $\bXi^N(\cdot)$ of the equation $\dot{\bXi}^N=\boldsymbol{F}(\bXi^N,\uu^N)$ with initial condition equal to 
$\bX$ has the following asymptotic 
behavior:
\begin{gather} \label{limit}
\lim_{t\to +\infty} x^N_{e_i}(t)=0 \qquad \lim_{t\to +\infty} y^N_{e_i}(t)=0\\
\lim_{t\to +\infty} z^N_{e_i}(t)=-1
\end{gather}
for $1\leq i \leq N$.
\end{proposition}

\proof
The proof relies on the fact that the restriction of $\bXi^N$ to the first $N$ components obeys to the dynamical system \eqref{eq: system}-\eqref{eq: control}. Then we can apply Theorem 
\ref{th: 1} and conclude that there exists an open dense subset  $A'$ of $\prod_{i=1}^N S^2$ such that for every $\bX\in \Sf$ with $\{X_{e_1},\ldots,X_{e_N}\}\in A'$ the 
solution $\bXi^N(\cdot)$ 
of the equation $\dot{\bXi}^N=\boldsymbol{F}(\bXi^N,\uu^N)$ with initial condition equal to $\bX$ satisfies the behavior described in 
\eqref{limit}, independently on the value of
$\{X_i\}_{k\geq N+1}$.
\eproof

\smallskip

Proposition \ref{prop: trunc} leads to the asymptotic pointwise convergence of the trajectories of $\bf(\bX,\uu)$ to $\bX^-$, where the notion of
``asymptotic pointwise'' convergence is weaker than the usual one, and is given in the following definition:

\begin{definition} \label{convergence}
The sequence of functions $\bXi^N(\cdot)$, with $\bXi^N(\cdot):\mathbb{R} \to \Sf$ for every $N$, 
converges asymptotically pointwise to the point $\bX\in \Sf$
if for every $\epsilon>0$ there exists an integer $\bar{N}>0$ such that for every $N\geq \bar{N}$ there exists a time $t=t(N,\epsilon)$ such that if $t\geq t(N,\epsilon)$ 
then $\boldsymbol{d}(\bXi^N(t),\bX)\leq \epsilon$.
\end{definition}

We can then state the following result.
 
\begin{theorem} \label{th: 3}
  There exists a  residual set $\mathcal{G}\subset \Sf$ such that for every $\bX\in \mathcal{G}$ there exists a sequence $\{\uu^N\}_N$ of controls 
  such that the sequence $\{\bXi^N\}_N$ of 
  solutions of the equation $\dot{\bXi}^N=\boldsymbol{F}(\bXi^N,\uu^N)$
with initial condition equal to $\bX$ converges asymptotically pointwise to $\bX^-$.
  \end{theorem}

\proof
Let $\bG^N\subset \prod_{i=1}^N S_i^2$ be the set of ``good initial conditions'' for the $N$-dimensional system, as defined in Theorem \ref{th: 1}, 
and let us define $\widehat{\bG}^N=\bG^N\times \prod_{i\geq N+1} S^2\subset \Sf$. Proposition \ref{prop: trunc}
states that the solution of the truncated system $\dot{\bXi}^N=\boldsymbol{F}(\bXi^N,\uu^N)$ with initial condition in $\widehat{\bG}^N$ has the limit \eqref{limit}. 
Since $\widehat{\bG}^N$ is an
open dense subset of $\Sf$ for every $N$, and
 $\Sf$ has the Baire property (\cite{kelley}), then $\mathcal{G}=\cap_{N}\widehat{\bG}^N$ is a dense subset of $\Sf$.

Let $\bX\in \mathcal{G}$ and fix $\epsilon>0$. For some integer $N$ such that $2^{-N+1}<\epsilon$, consider the truncated feedback  $\uu^N$, defined as above, and the corresponding
trajectory $\bXi^N$ with $\bXi^N(0)=\bX$. Since $\bX\in \widehat{\bG}^N$, by Proposition~\ref{prop: trunc}
there exists a time
$t=t(N,\epsilon)$ such that for $t\geq t(N,\epsilon)$  
it holds $\sum_{i=1}^N 2^{-i} |X^N_{e_i}(t)-(0,0,-1)^T|\leq \epsilon/2$, then, since $\sum_{k\geq N+1}2^{-i}|X^N_{e_i}(t)-(0,0,-1)^T|\leq 2^{-N} \leq\epsilon/2$, we get that 
$\boldsymbol{d}(\bXi^N(t),\bX^-)\leq \epsilon$ for $t\geq t(N,\epsilon)$.
\eproof

\section{Closed-loop simulations}\label{sec: simu}

Let $\mathcal{E}$ be a collection of $N=30$ randomly chosen points contained in the interval $[1,4]$, and we consider $N$ randomly chosen initial 
conditions $X_e^0$ with $z_e^0\in[0.8,1]$ and 
$|X_e^0|=1$.
We perform closed-loop simulation of the dynamical system \eqref{eq: system} with feedback control $u_1(t)=\sum_{i=1}^N y_{e_i}(t)$ and 
$u_2(t)=\sum_{i=1}^N x_{e_i}(t)$, up to a final time $T=20000$.

\begin{figure}[h!]
        \begin{subfigure}[b]{0.30\textwidth}
                \includegraphics[width=\linewidth]{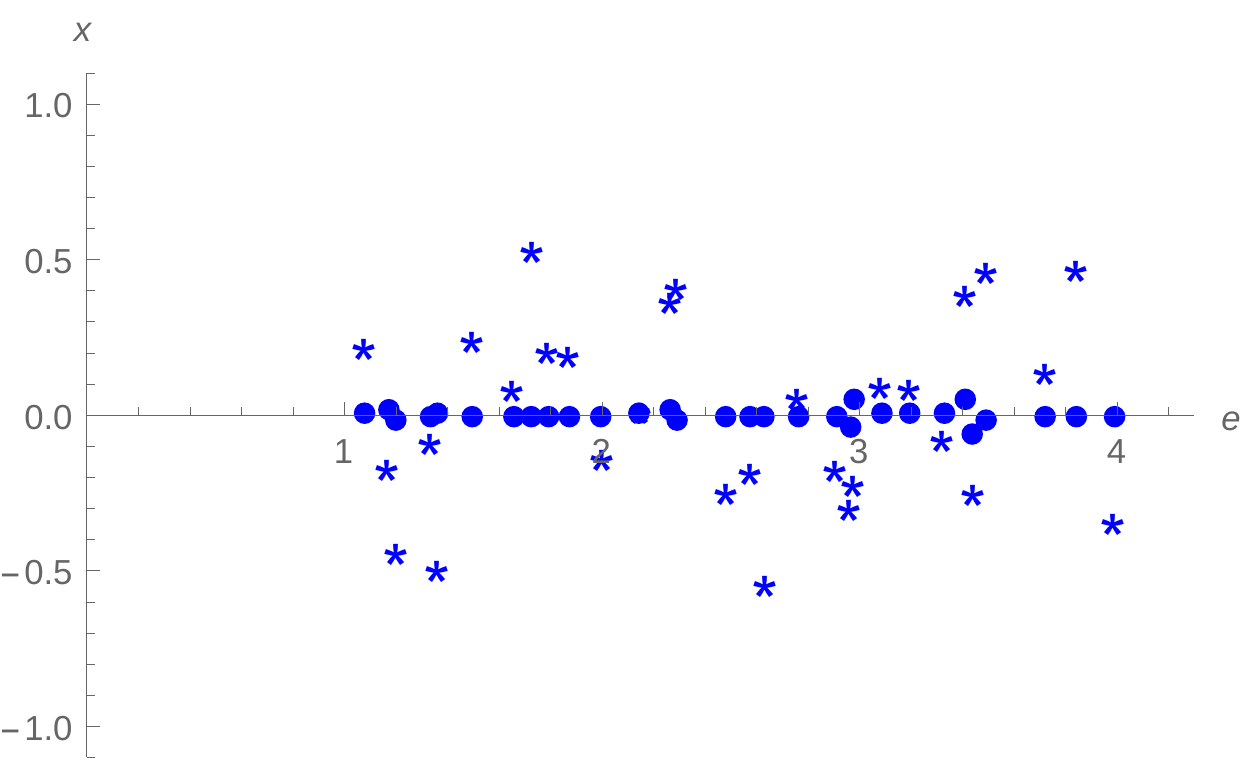}
                \caption{$x_e$}
                \label{fig: X}
        \end{subfigure}%
        \begin{subfigure}[b]{0.30\textwidth}
                \includegraphics[width=\linewidth]{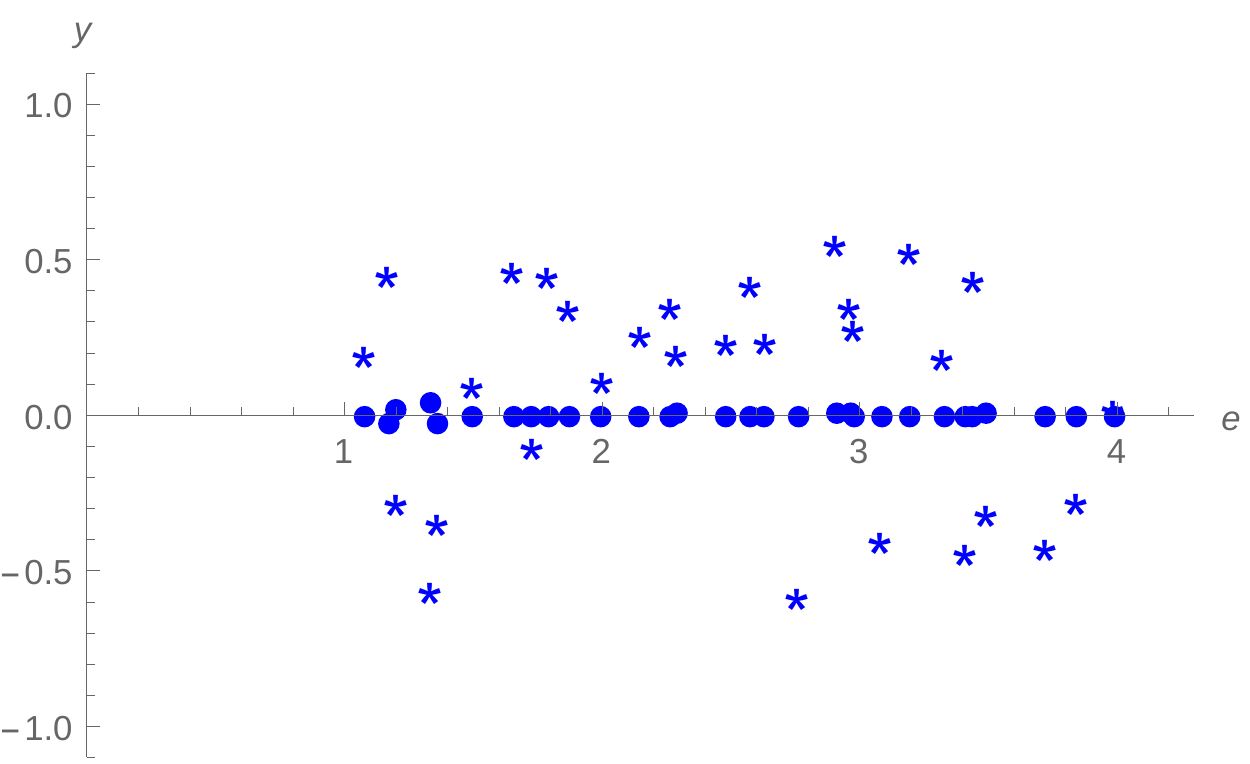}
                \caption{$y_e$}
                \label{fig: Y}
        \end{subfigure}%
        \begin{subfigure}[b]{0.30\textwidth}
                \includegraphics[width=\linewidth]{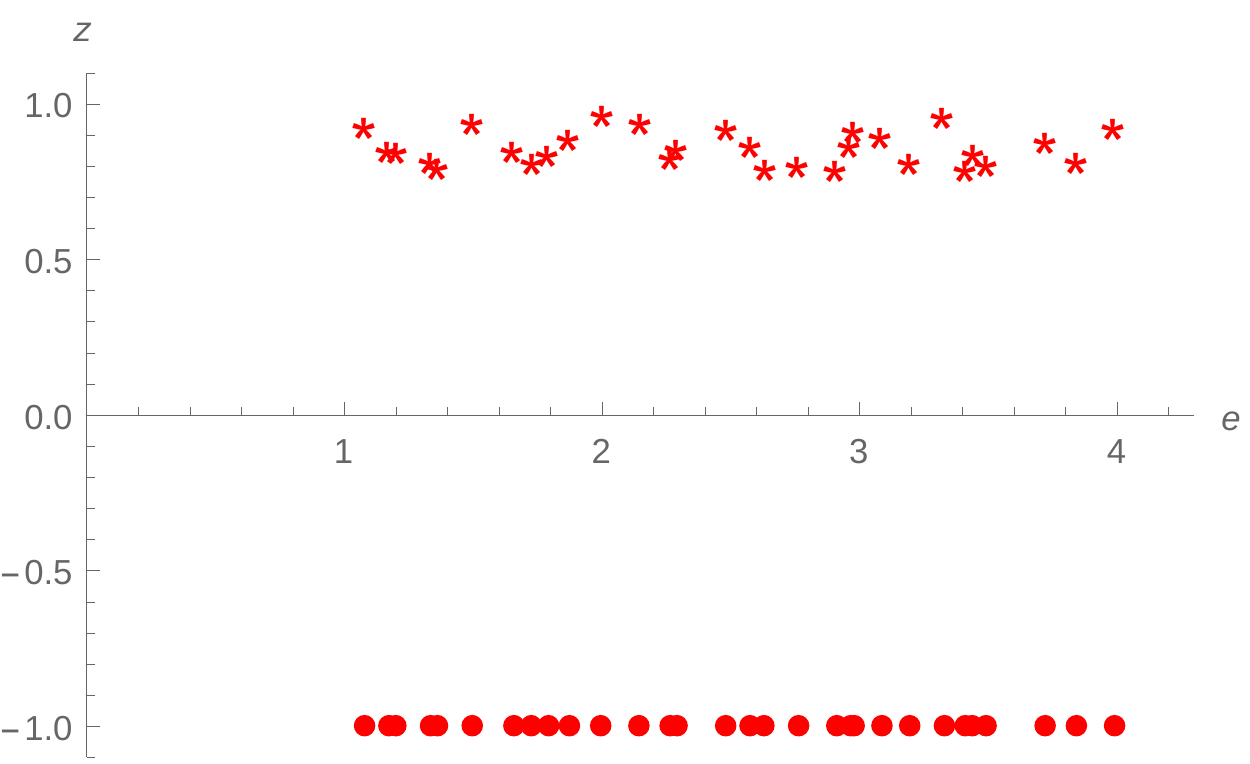}
                \caption{$z_e$}
                \label{fig: Z}
        \end{subfigure}
        \caption{Initial and final states with respect to different values of frequencies (the stars denote the initial point, the bullets the final point).}\label{fig: XYZ}
\end{figure}

In Figure \ref{fig: XYZ} we show the convergence to the target point of the collections $X_e$, $e\in \mathcal{E}$. 

\begin{figure}[h!]
        \begin{subfigure}[b]{0.45\textwidth}
                \includegraphics[width=\linewidth]{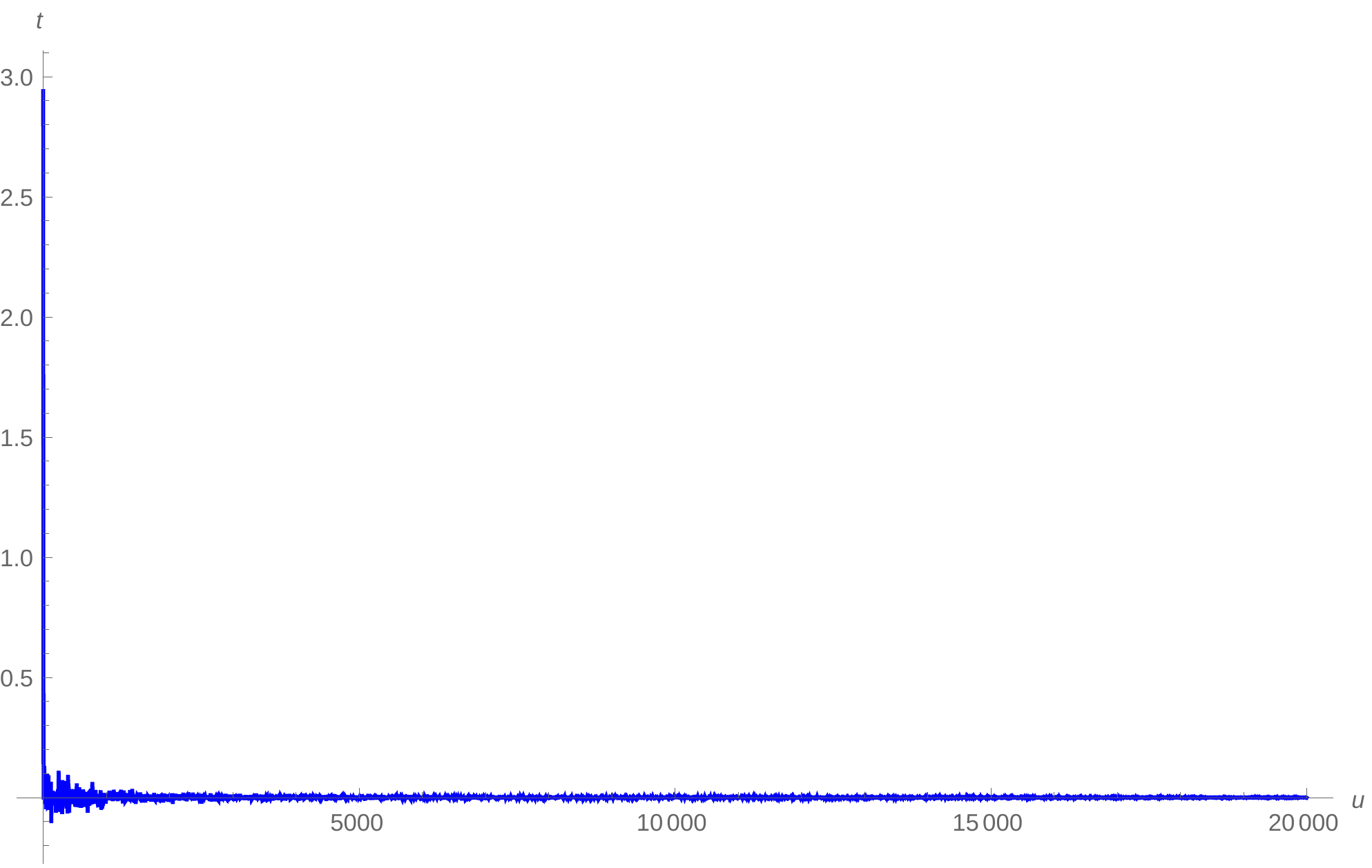}
                \caption{$u_1$}
                \label{fig: u1}
        \end{subfigure}%
        \begin{subfigure}[b]{0.45\textwidth}
                \includegraphics[width=\linewidth]{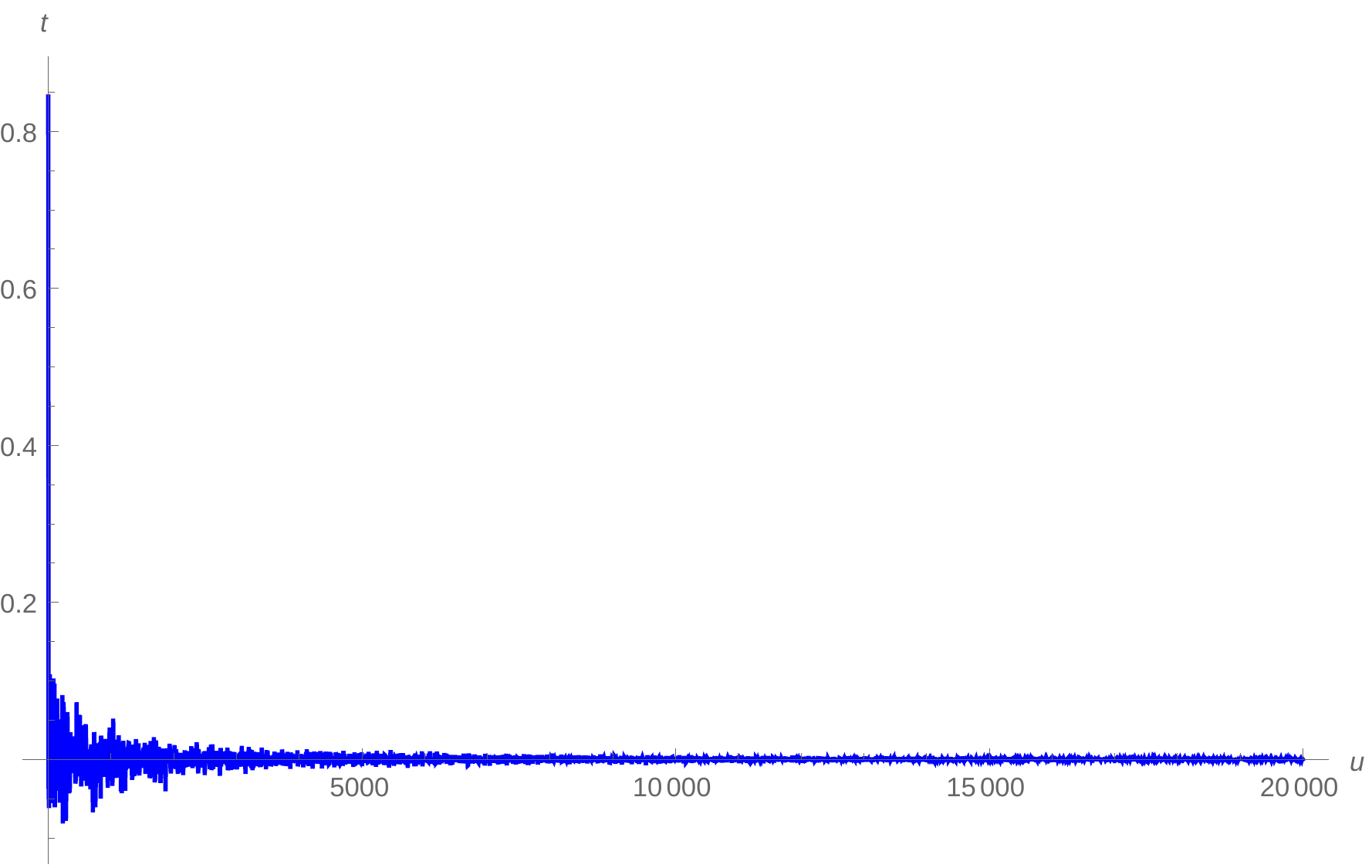}
                \caption{$u_2$}
                \label{fig: u2}
        \end{subfigure}
        \caption{Time evolution of the control function}\label{fig: U}
\end{figure}

\begin{figure}[h!]
\begin{subfigure}[b]{0.45\textwidth}
                \includegraphics[width=\linewidth]{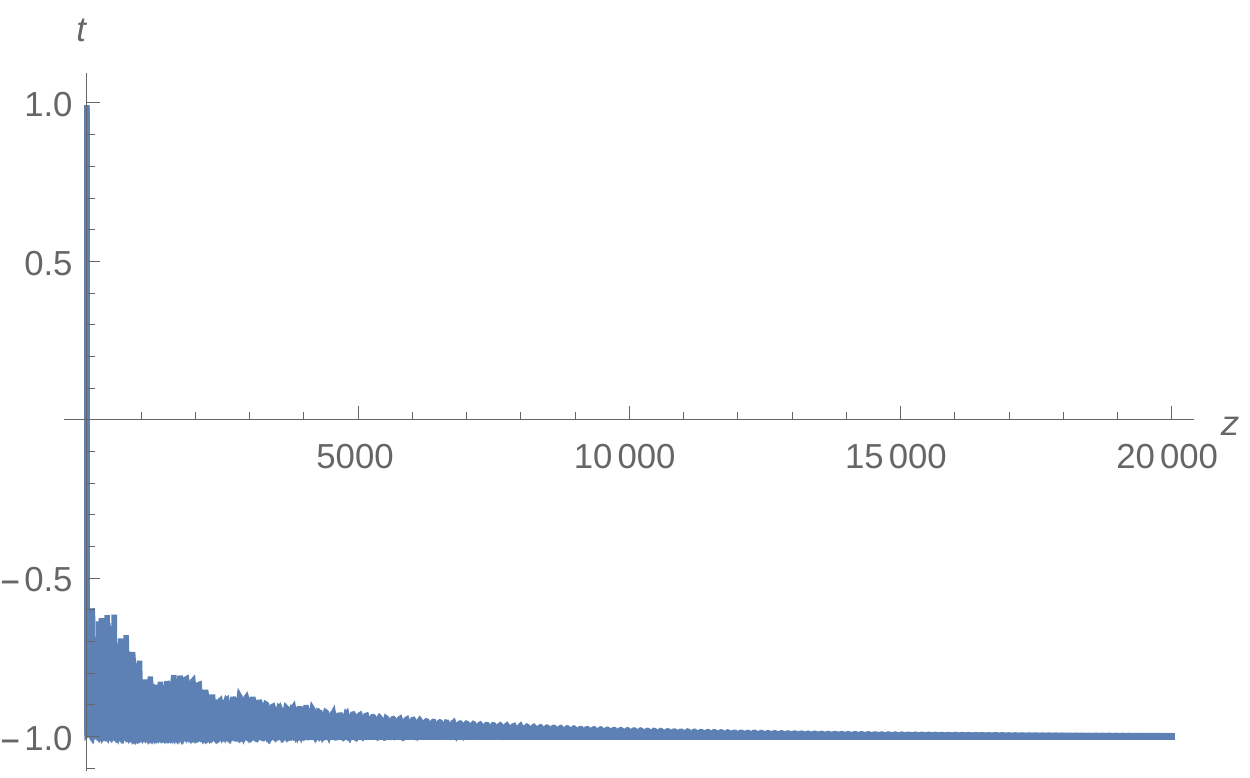}
                \caption{$z_e(t)$ for different values of $e$}
                \label{fig: Zdecay}
                \end{subfigure}%
                \begin{subfigure}[b]{0.45\textwidth}
                \includegraphics[width=\linewidth]{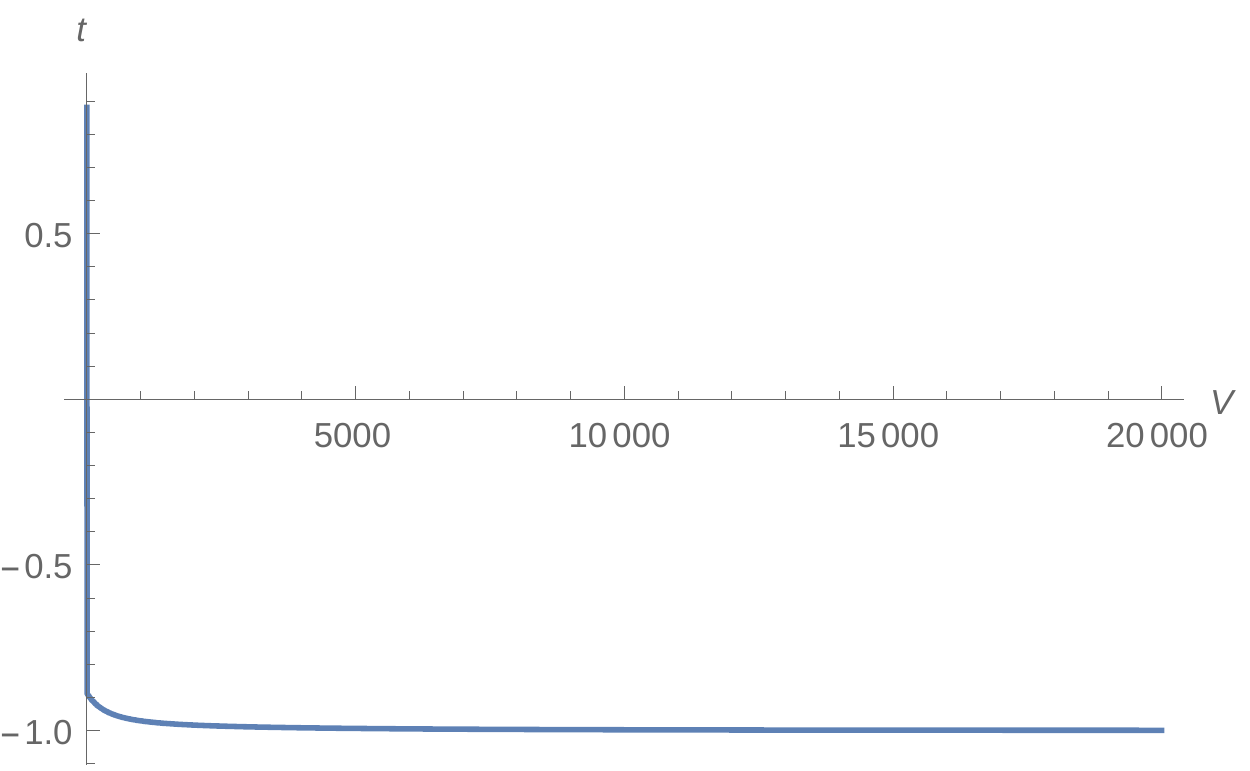}
                \caption{Lyapunov function (normalized)}
                \label{fig: lyap}
        \end{subfigure}\caption{Time evolution of $z_e$ and the Lyapunov function}
\end{figure}

Figure \ref{fig: U} plots the time evolution of the feedback control function, while in Figures \ref{fig: Zdecay} and \ref{fig: lyap} 
we plot respectively the values of the last coordinate $z_e(t)$, for all $e\in \mathcal{E}$,
and of the Lyapunov function $V(t)$, normalized by $N$. 

\begin{figure}[h!]
\begin{subfigure}[b]{0.45\textwidth}
                \includegraphics[width=\linewidth]{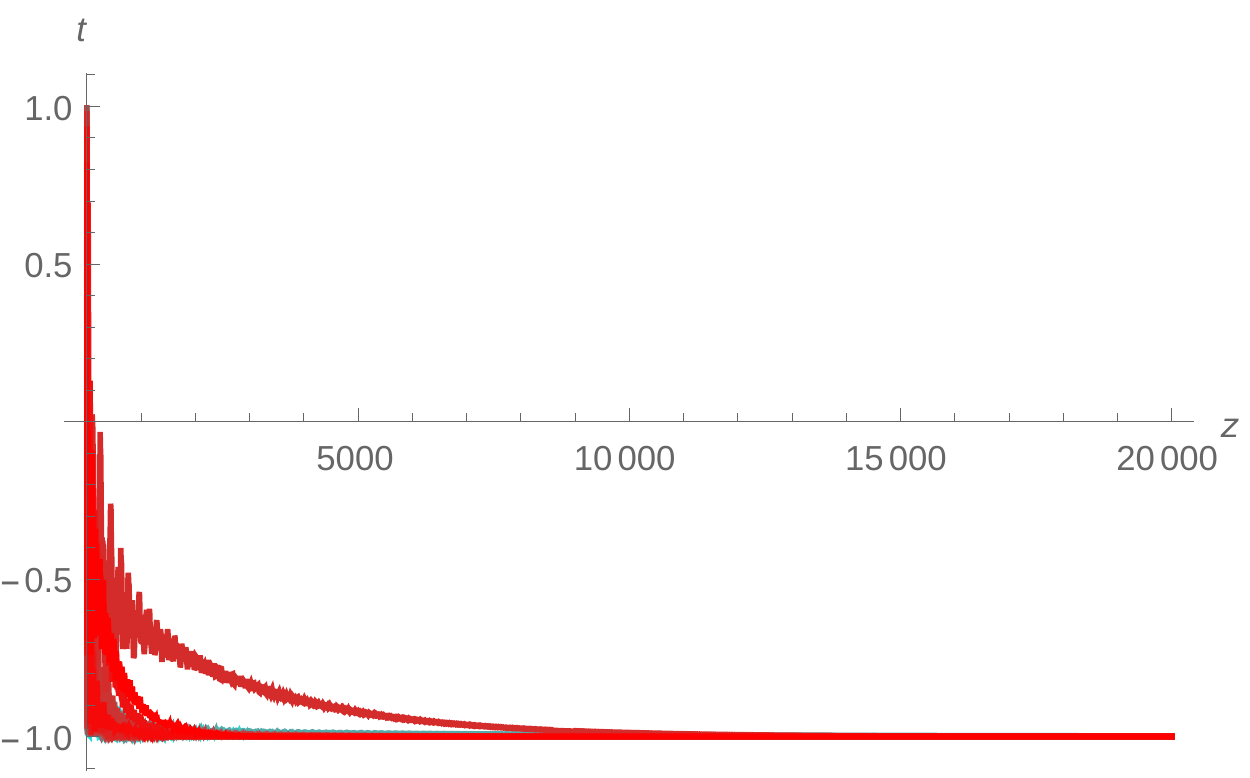}
                \caption{$z_e(t)$ for different values of $e$}
                \label{fig: Zdecayw}
                \end{subfigure}%
                \begin{subfigure}[b]{0.45\textwidth}
                \includegraphics[width=\linewidth]{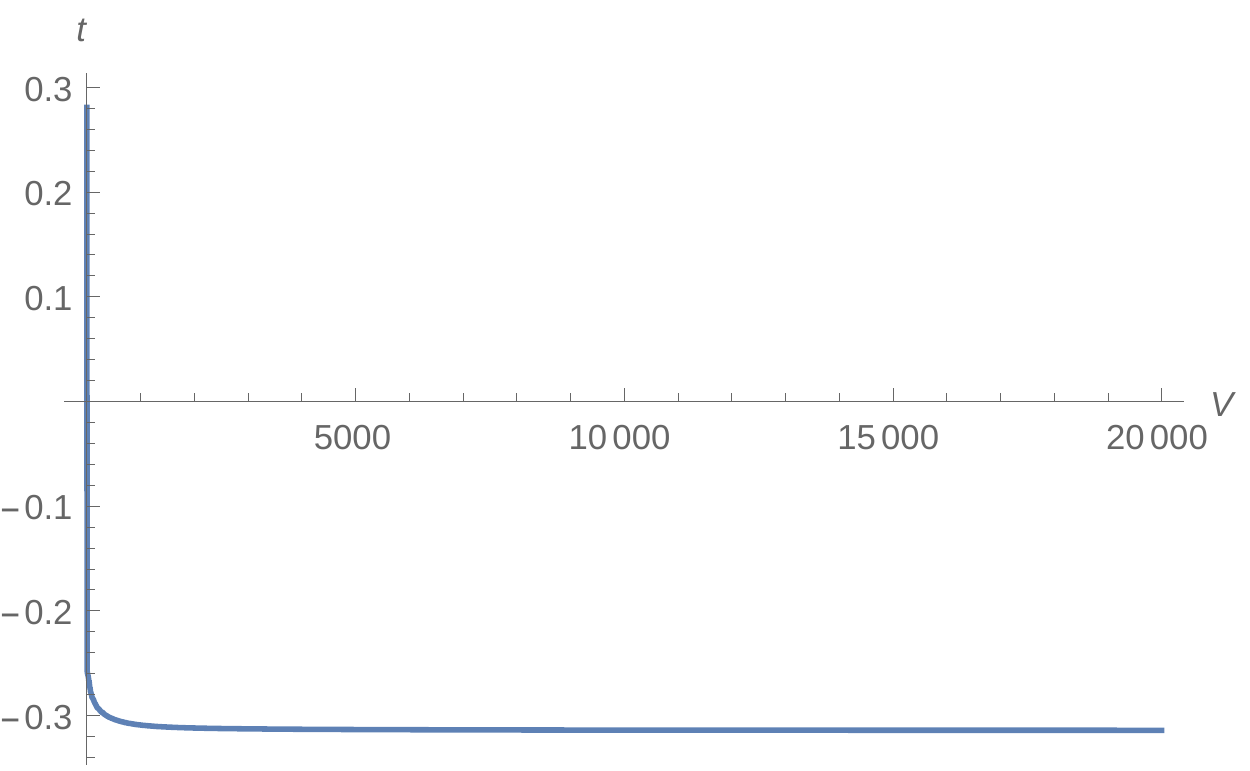}
                \caption{Lyapunov function (normalized)}
                \label{fig: lyapw}
        \end{subfigure}\caption{Time evolution of $z_e$ and the Lyapunov function}
\end{figure}

\medskip
We then take the same collection $\mathcal{E}$ as before, and we consider $N$ randomly chosen initial conditions $X_e^0$ with $z_e^0\in[0.8,1]$ and 
$|X_e^0|=1$.
We now perform closed-loop simulation of the dynamical system \eqref{eq: system} with feedback control $u_1(t)=\sum_{i=1}^N w_i y_{e_i}(t)$ and 
$u_2(t)=\sum_{i=1}^N w_i x_{e_i}(t)$, with $w_i=(1.1)^{-i}$, up to a final time $T=20000$.
The purpose of this new run is to visualize the influence of the weights $w_i$ on the convergence of the systems. As we can see from Figure~\ref{fig: Zdecayw},
the weights slow down the convergence of the systems (this cannot be seen from Figure~\ref{fig: lyapw}, since the slower components in the Lyapunov function
are multiplied by a small weight).

\begin{figure}[h!]
        \begin{subfigure}[b]{0.30\textwidth}
                \includegraphics[width=\linewidth]{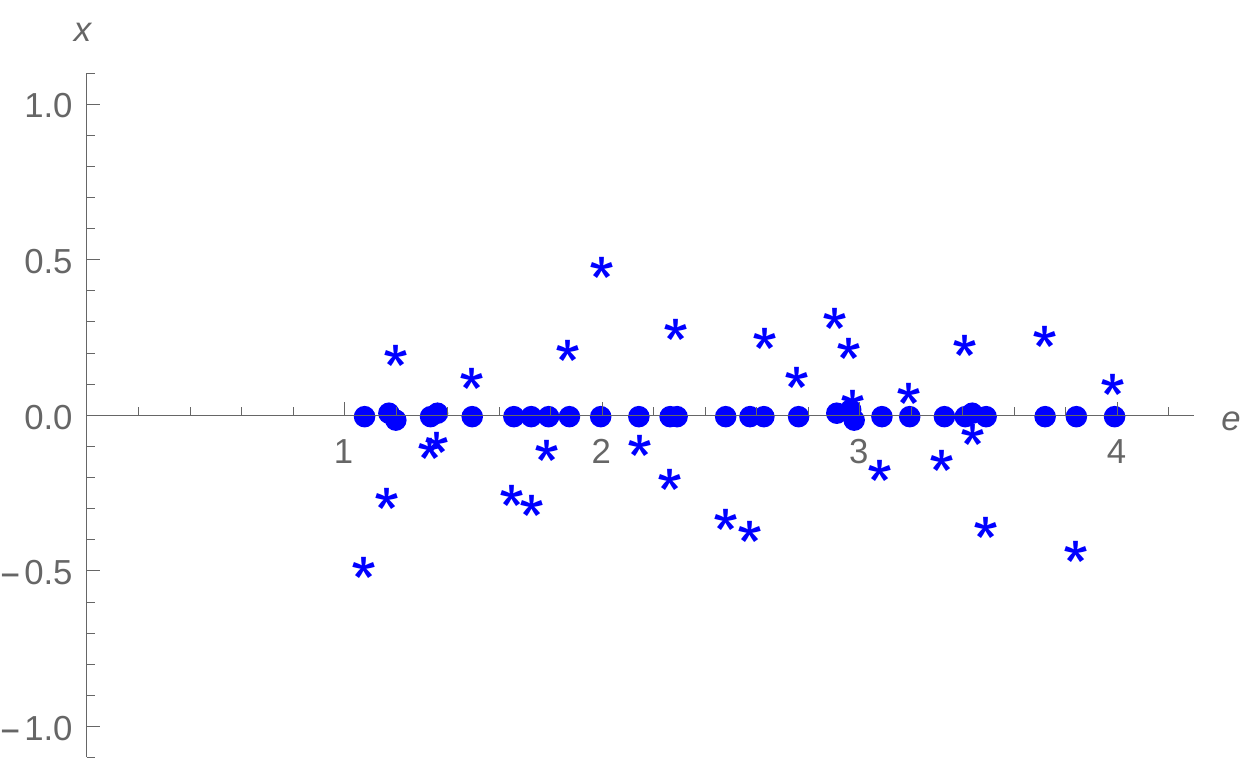}
                \caption{$x_e$}
                \label{fig: Xw}
        \end{subfigure}%
        \begin{subfigure}[b]{0.30\textwidth}
                \includegraphics[width=\linewidth]{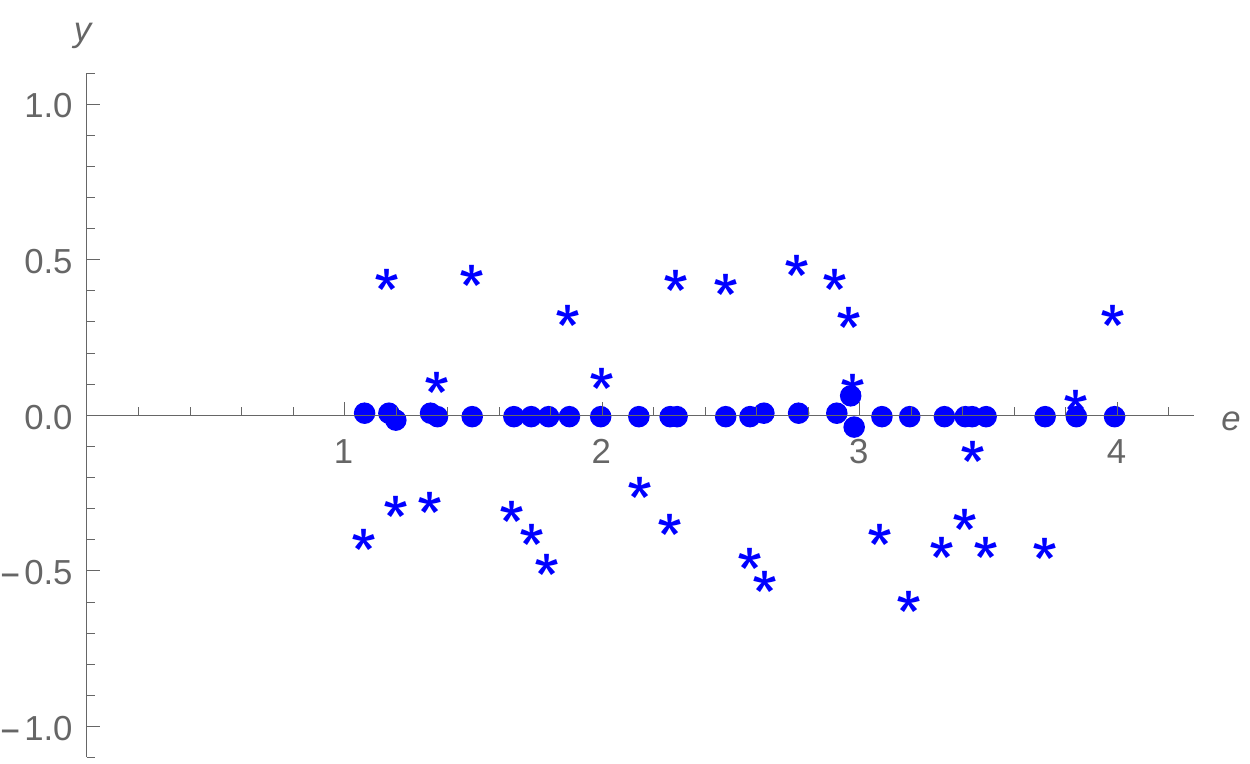}
                \caption{$y_e$}
                \label{fig: Yw}
        \end{subfigure}%
        \begin{subfigure}[b]{0.30\textwidth}
                \includegraphics[width=\linewidth]{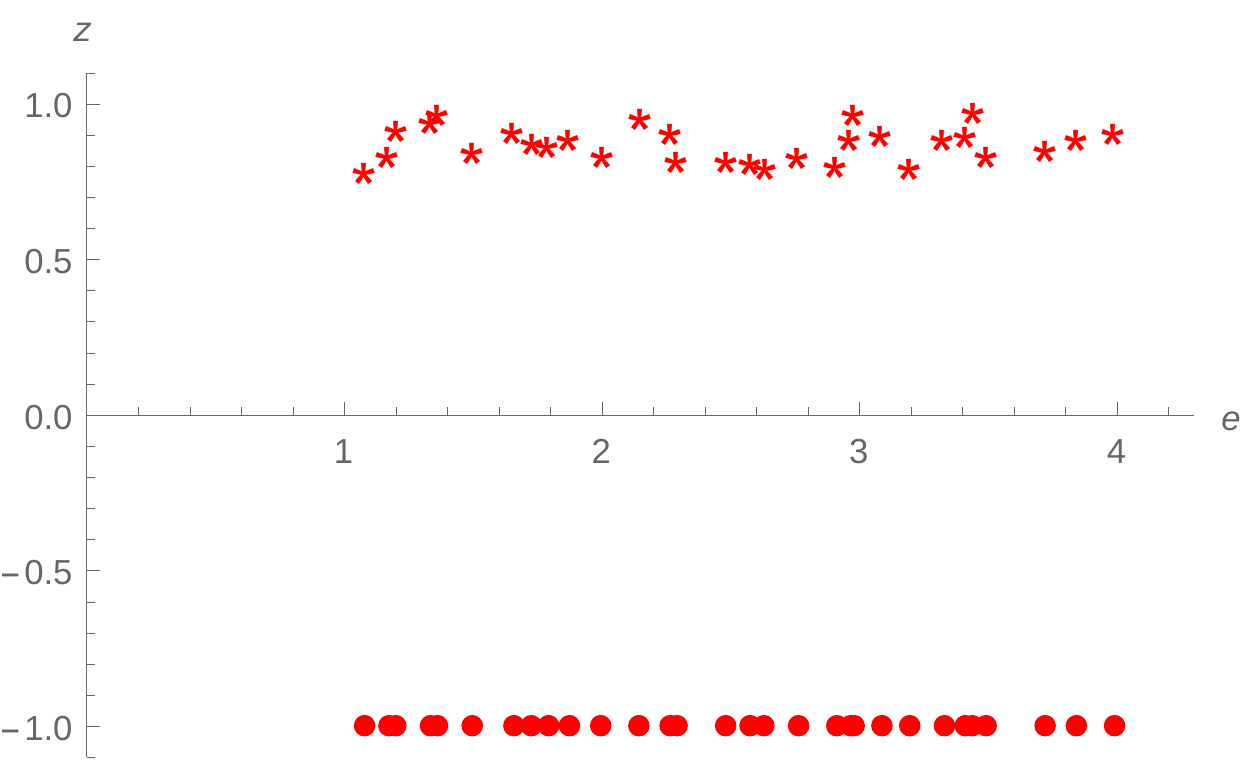}
                \caption{$z_e$}
                \label{fig: Zw}
        \end{subfigure}
        \caption{Initial and final states with respect to different values of frequencies (the stars denote the initial point, the bullets the final point).}
        \label{fig: XwYwZw}
\end{figure}

In Figure \ref{fig: XwYwZw} we show the convergence to the target point of the collections $X_e$, $e\in \mathcal{E}$. 

\begin{figure}[h!]
        \begin{subfigure}[b]{0.45\textwidth}
                \includegraphics[width=\linewidth]{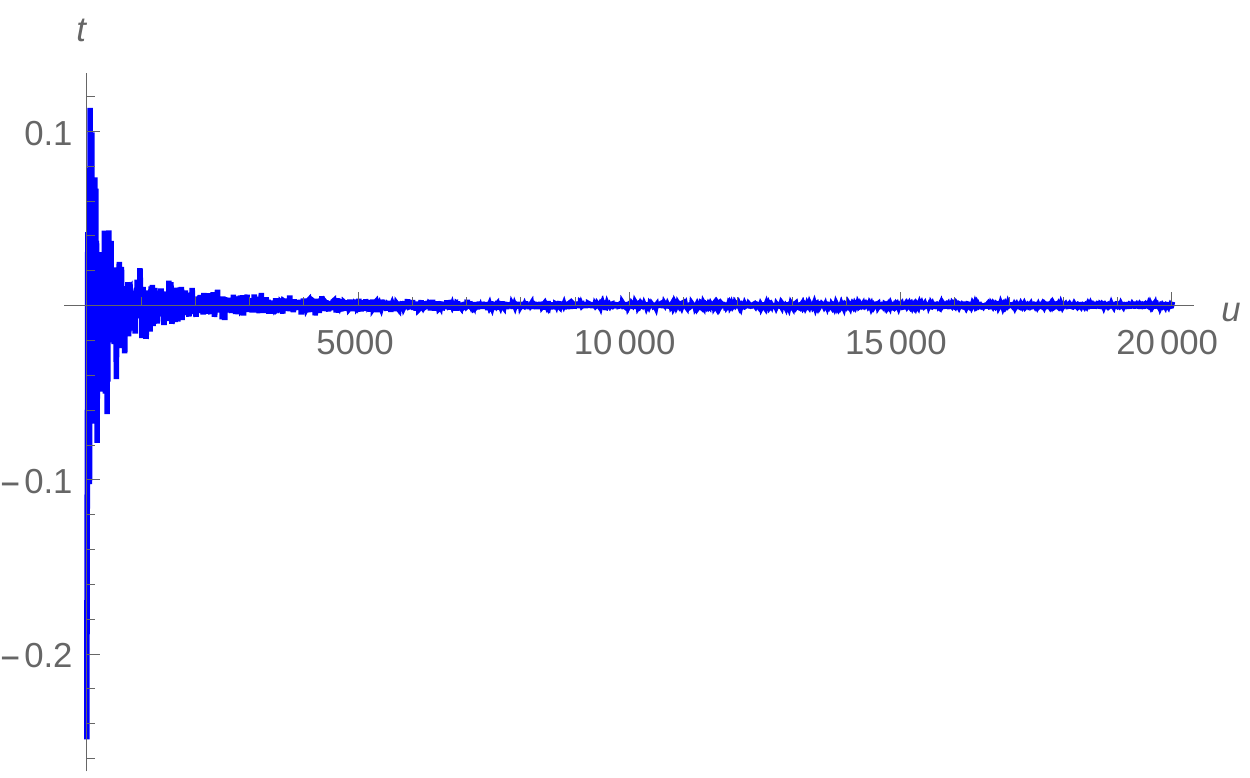}
                \caption{$u_1$}
                \label{fig: u1w}
        \end{subfigure}%
        \begin{subfigure}[b]{0.45\textwidth}
                \includegraphics[width=\linewidth]{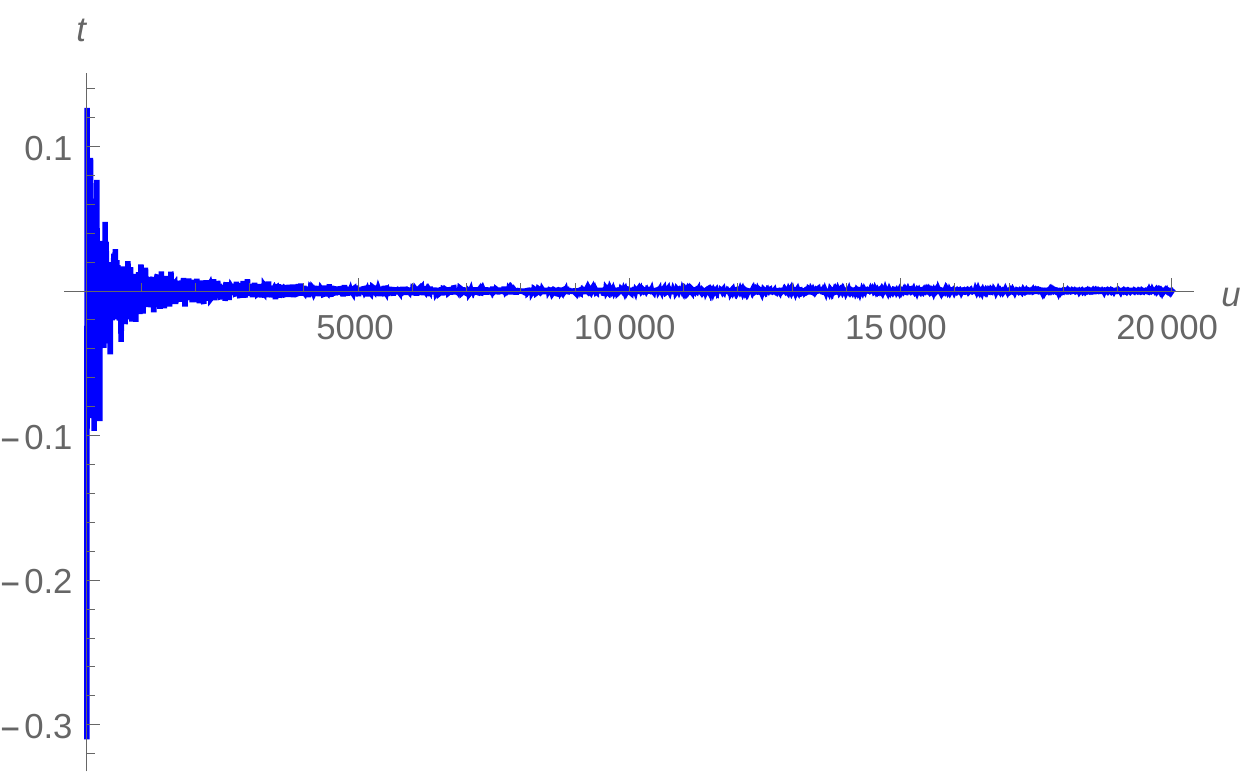}
                \caption{$u_2$}
                \label{fig: u2w}
        \end{subfigure}
        \caption{Time evolution of the control function}\label{fig: Uw}
\end{figure}

As above, in Figure \ref{fig: Uw} we plot the time evolution of the feedback control function, in Figures \ref{fig: Zdecayw} $z_e(t)$, 
and in \ref{fig: lyapw} 
 the Lyapunov function $V(t)$, normalized by $N$. 

\section{Conclusions}

In this paper, we have investigated the stabilization of an ensemble of non-interacting half-spins to the uniform state $-1/2$ (represented by the state $(0,0,-1)$ in the Bloch sphere); in particular,
we provided a feedback control that stabilizes a generic initial condition to the target state, asymptotically in time. 

In the finite-dimensional case, we remark a close link between the proposed control \eqref{eq: control} and the \emph{radiation damping effect (RDE)} (see for instance 
\cite{bloem,augustine} for a detailed description of the phenomenon). 
In an NMR setup, the radiation damping is a reciprocal interaction between the spins and the radio-frequency source (a coil):
this coupling can be taken into account by adding
a non-linear term to the uncontrolled
Bloch equation (see for instance \cite{altafinicory,zhang} and references therein). 
In particular, in our notations the uncontrolled Bloch equation with RDE reads 

\begin{equation} \label{eq: system wd}
\begin{cases}
\dot{x}_e =  -e_iy_{e_i} - \ell z_{e_i} \overline{X} \\ 
\dot{y}_e=
e x_{e_i} -\ell z_{e_i}\overline{Y}\\
\dot{z}_e
=\ell \big( \overline{X}^2 +\overline{Y}^2 \big)
\end{cases},
\end{equation}
where $\ell$ is the radiation damping rate (depending on the apparatus) and $\overline{X}=\frac{1}{p}\sum_{i=1}^p x_{e_i},\ \overline{Y}=\frac{1}{p}\sum_{i=1}^p y_{e_i}$ are the average 
values of the magnetization. 
The analysis carried out in Section~\ref{sec: finite} applies also in this case, with the only difference that $\bX^-$ is a repeller and  $\bX^+$ is an attractor of
equation~\eqref{eq: system wd}. This gives a rigorous justification of the stabilizing properties of RDE.

If we want to exploit RDE for stabilizing the system towards $\bX^-$, it is sufficient to
invert the $z$-component of the magnetic field: this yields a change of the sign of  the right-hand side of equation~\eqref{eq: system wd}, 
thus, up to a change in the sign of the frequencies
(which does not affect the dynamics, being the set $\mathcal{E}$ arbitrary) and to a multiplicative factor $\ell/p$ on the control, 
we obtain the dynamical system \eqref{eq: system}-\eqref{eq: control}. 
The multiplicative factor $\ell/p$
affects only the magnitude of the real part of the eigenvalues (see equation \eqref{eigenvalue fin}), that is, the rate of convergence towards the equilibria.

If it is not possible to invert the $z$-component of the magnetic field, so that the RDE tends to stabilize the system to  $\bX^+$, the stabilization to  $\bX^-$ can be still achieved by 
choosing a sufficiently strong control (see for instance \cite{altafinicory} for a similar result in the single spin case).

%

\medskip
In the countable case, we use the same approach to provide a sequence of (continuous bounded) feedback controls which asymptotically stabilizes, 
according to the notion of convergence given in Definition~\ref{convergence}, 
a generic set of initial conditions.

Concerning the case where $\mathcal{E}$ is an interval, and $\bX\in L^2(\mathcal{E},S^2)$, the question addressed in \cite{beauchard-coron-rouchon2010} about 
controllability of the system by means of bounded controls is still left open. 
This topic makes the subject of further investigations of the authors.

\bibliographystyle{plain}
\bibliography{biblio-ens}

\end{document}